\documentclass[11pt,reqno]{amsart}
\usepackage[T1]{fontenc}
\usepackage{array}
\usepackage{lmodern,amsmath,amssymb,mathtools,microtype,enumitem}
\usepackage[letterpaper,margin=1.04in,top=0.96in,bottom=0.96in]{geometry}
\usepackage[hidelinks,bookmarksnumbered=true,pdfauthor={Hengzhi He and Guang Cheng},pdftitle={The smallest singular value of a random sign matrix: the Spielman--Teng conjecture}]{hyperref}
\allowdisplaybreaks[2]
\newtheorem{theorem}{Theorem}[section]
\newtheorem{proposition}[theorem]{Proposition}
\newtheorem{lemma}[theorem]{Lemma}
\newtheorem{corollary}[theorem]{Corollary}
\theoremstyle{definition}
\newtheorem{definition}[theorem]{Definition}
\newtheorem{remark}[theorem]{Remark}
\numberwithin{equation}{section}
\newcommand{\R}{\mathbb R}
\newcommand{\Z}{\mathbb Z}
\newcommand{\E}{\mathbb E}
\newcommand{\PP}{\mathbb P}
\newcommand{\ind}{\mathbf 1}
\newcommand{\Acal}{\mathcal A}
\newcommand{\Gcal}{\mathcal G}
\newcommand{\Hcal}{\mathcal H}
\newcommand{\smin}{s_{\min}}
\newcommand{\tr}{\operatorname{tr}}
\newcommand{\diag}{\operatorname{diag}}
\newcommand{\dist}{\operatorname{dist}}
\newcommand{\rank}{\operatorname{rank}}
\newcommand{\vol}{\operatorname{vol}}
\newcommand{\sph}{\mathrm{sph}}

\newcommand{\one}[1]{\ind_{\{#1\}}}
\newcommand{\step}[2][.]{\smallskip\noindent\emph{#2#1}\enspace}
\newcommand{\ga}{\phi_0}
\newcommand{\sphere}{S^{n-1}}
\newcommand{\cflat}{c_{\flat}}
\setlist[enumerate]{itemsep=1pt,topsep=3pt,parsep=0pt}
\setlist[itemize]{itemsep=1pt,topsep=3pt,parsep=0pt}
\title[The Spielman--Teng conjecture for random sign matrices]{The smallest singular value of a random sign matrix:\\ the Spielman--Teng conjecture}
\author{Hengzhi He}
\author{Guang Cheng}
\begin{document}
\begin{abstract}
Let \(A\) be an \(n\times n\) matrix with independent uniform entries in \(\{-1,1\}\). We prove that
\[
 \PP\bigl(\smin(A)\le \varepsilon n^{-1/2}\bigr)\le \varepsilon+e^{-cn}
 \qquad(\varepsilon\ge 0)
\]
for an absolute \(c>0\), as conjectured by Spielman and Teng. For \(e^{-cn}\le\varepsilon\le n^{-\theta}\), with \(\theta>0\) absolute and \(n\) sufficiently large, the upper bound improves to
\(\varepsilon(1-1/(2n)+O(n^{-1.02}))\).
The proof does not give numerical values for the dimension threshold or the exponential rate. We obtain the first-order improvement by comparing a Gaussian surface integral for the sign ensemble with its counterpart for independent spherical rows. An arithmetic small-ball estimate transfers this comparison to exponentially small scales. At polynomial scales, weighting by a coordinate of the least left singular vector yields a negative quadratic term in the crossing density.
\end{abstract}
\maketitle

\section{Introduction}
\subsection{The result}
Let \(A\) be an \(n\times n\) matrix whose entries are independent uniform signs, and write
\(\smin(A)=\min_{|x|=1}|Ax|\). If \(G\) has independent standard real Gaussian entries, then
\(\PP(\smin(G)\le t/\sqrt n)\le t\) for every \(t\ge0\)
\cite{edel,szarek}. Remark~\ref{rm:Gaussian} gives a short derivation.
In their 2002 ICM survey, Spielman and Teng \cite[Conjecture~2, p.~605]{st} proposed the same inequality for a sign matrix, with an exponentially small allowance for singularity. They formulated it in terms of the inverse norm, with \(\|A^{-1}\|=\infty\) when \(A\) is singular. Taking limits in the threshold gives the equivalent formulation below.
\begin{theorem}\label{th:main}
There is an absolute constant \(c>0\) such that for all \(n\ge1\) and \(t\ge0\),
\[
 \PP\bigl(\smin(A)\le t/\sqrt n\bigr)\le t+e^{-cn}.
\]
\end{theorem}
Define \(X_n=\sqrt n\,\smin(A)\) and \(F_n(t)=\PP(X_n\le t)\).
\begin{theorem}\label{th:sharper}
There are absolute constants \(c,\theta>0\) and \(C,n_0<\infty\) such that for \(n\ge n_0\) and \(e^{-cn}\le t\le n^{-\theta}\),
\[
 F_n(t)\le t\bigl(1-\tfrac1{2n}+Cn^{-1.02}\bigr).
\]
\end{theorem}

There are two quantifiers worth keeping in view. Theorem~\ref{th:main} covers every dimension, including at \(t=0\), where its left side is the singularity probability. The coefficient in Theorem~\ref{th:sharper} is a one-sided bound for sufficiently large \(n\). Our estimates require \(n^{\theta/2000}\) to dominate \(\log n\), as well as a qualitative convergence result for singular values. With the exponents used here the required dimensions can be extremely large. Section~\ref{sec:finish} handles the finitely many smaller dimensions and explains how that cutoff enters \(c\).

\subsection{Previous work}
Rudelson \cite{rud} and Tao and Vu \cite{taoinv} obtained early polynomial bounds for discrete matrices. Rudelson and Vershynin proved an upper bound \(Ct+e^{-cn}\) for the square problem \cite{rvsmall} and established its rectangular counterpart \cite{rvrect}. Tao and Vu \cite{taodist} proved universality for \(X_n\) with polynomial error. Their theorem gives the conjectured upper bound for \(t\ge n^{-c'}\) with a sufficiently small absolute \(c'>0\). 

For iid, centered, variance-one subgaussian entries, Sah, Sahasrabudhe and Sawhney \cite{sss} obtained the upper bound \((1+o(1))t+e^{-\Omega(n)}\) uniformly in \(t\). They also pointed to expansions for complex ensembles \cite{egp} suggesting that an exact coefficient of one could require distribution-specific information. 

For the Rademacher law, Hu, Song and Wu \cite{hsw} obtain \(\PP(A\text{ is singular})=(2+o(1))n^22^{-n}\) and a positive-threshold bound with linear term \(Ct\). The constant one at positive thresholds is the subject of Theorem~\ref{th:main}. Our comparison uses both the fixed length of each row and the characteristic function of its projections.

\subsection{Why the proof uses three scales}
Remove the last row of \(A\) and call the remaining matrix \(B\). When \(B\) has full row rank, it has a single missing direction, its unit normal \(v\). If the last row \(X^\top\) lies in the row space of \(B\), then \(A\) is singular. Its signed distance from this row space is \(\ell=\langle v,X\rangle\). For a positive threshold, the minimizing direction in \(\min_{|x|=1}|Ax|\) can move away from \(v\). A Schur complement gives the following exact criterion for a near-singular matrix.
\begin{equation}\label{eq:outline}
 X_n\le t \quad\Longleftrightarrow\quad |\ell|\le t W_t,
 \qquad 0\le t<\sqrt n\,s_1(B).
\end{equation}
The random width \(W_t\) depends on the tangential projections of the same last row. It can therefore be correlated with \(\ell\). Also, conditional on \(B\), \(\ell\) takes only finitely many values. These features explain why a density approximation for \(\ell\) alone would leave the leading coefficient undetermined.

\step[?]{Where does the improvement over one come from}
At the reference scale \(t_0=n^{-25}\), we replace a sharp probability cutoff by a Gaussian average of \(X_n\). The image of the unit sphere under \(A\) has a short semiaxis precisely when \(A\) has a small singular value. We bound the Gaussian average by integrating over the part of this image near zero. In each direction its integrand factors into the \(n\)-th power of a one-row projection density and an expected tangential volume. Comparing the volume with independent spherical rows gives the factor \(1-5/(4n)+O(n^{-2})\). The angular variation of the projection density costs at most \(1+3/(4n)+O(n^{-3/2})\). Multiplying the two leaves a margin of \(1/(2n)\) below one.

\step[?]{How can a bound on an average control individual thresholds}
For \(e^{-\cflat n}\le t\le t_*=n^{-20}\), we prove that \(F_n(t)/t\) differs from a single number \(a_n\) by at most \(Cn^{-1.05}\). Conditionally on \(B\), the event in \eqref{eq:outline} is a thin interval whose width is proportional to \(t\). After smoothing the normal and tangential coordinates together, the mass of that interval is its width times the joint density at normal coordinate zero, up to a controlled error. Since \(t_0\) lies inside this range, the Gaussian average bounds \(a_n\). This argument also allows an atom at zero.

\step[?]{Why change the argument when \(t\) is larger}
At polynomial thresholds a first-order approximation to the interval mass is too coarse to preserve a margin of order \(t/n\). Let \(u\) be a least left singular vector chosen as in Lemma~\ref{le:exchange}. Row symmetry gives the exact identity \(F_n(t)=\E[nu_n^2\one{X_n\le t}]\). The weight cancels part of the Jacobian when we change from the normal coordinate to the threshold at which the event first occurs. The resulting density has a negative quadratic correction on average. Integrating it gives enough room to include the blocks with an unusually small singular value. At \(t_*\), the weighted and unweighted estimates concern the same \(F_n(t_*)\), so their bounds can be matched there.

\clearpage
\subsection{A guide to the proof}\label{sec:guide}
\noindent The proof of the sharp coefficient can be read by following the four estimates below. The two reference scales play different roles. We use \(t_0=n^{-25}\) to compute the coefficient and \(t_*=n^{-20}\) to match the two probability estimates.

\smallskip
\begingroup\normalsize\renewcommand{\arraystretch}{1.4}
\noindent\begin{tabular}{@{}p{.23\linewidth}@{\hspace{.035\linewidth}}>{\raggedright\arraybackslash}p{.71\linewidth}@{}}
\hline
\textbf{Estimate} & \textbf{Inputs and conclusion}\\\hline
Gaussian average at \(t_0\) &
Lemmas~\ref{le:area} and \ref{le:zonoid} compare the surface area with a spherical model. Appendix~\ref{app:sph} computes its coefficient, and Lemma~\ref{le:angle} bounds the angular moment. Together they prove Proposition~\ref{pr:area}.\\
Flatness up to \(t_*\) &
Corollary~\ref{co:slab} and Lemma~\ref{le:rare} control the last row on a thin slab. Proposition~\ref{pr:flat} gives a common slope. Corollary~\ref{co:slope} uses the Gaussian average to bound it.\\
Extension from \(t_*\) to \(n^{-\theta}\) &
Lemmas~\ref{le:smallblock}--\ref{le:compare} express the weighted probability through crossings. Lemmas~\ref{le:high}--\ref{le:curv} control the Fourier tail and prove the sign of the quadratic correction. Corollary~\ref{co:poly} matches the slopes at \(t_*\).\\
Remaining thresholds &
Below \(e^{-\cflat n}\) use monotonicity. Above \(n^{-\theta}\), Tao--Vu universality and the explicit limiting distribution give the required bound in Section~\ref{sec:finish}.\\\hline
\end{tabular}
\endgroup

\medskip\noindent\textbf{How the parameters are chosen.}
The universality error fixes \(\theta=\min(c_{\rm TV}/4,1/100)\). It in turn fixes \(\kappa=\theta/20\) and \(\tau=\theta/1000\). The remaining powers serve the following estimates.
\par\smallskip
\begingroup\normalsize\renewcommand{\arraystretch}{1.30}
\noindent\begin{tabular}{@{}p{.31\linewidth}@{\hspace{.035\linewidth}}>{\raggedright\arraybackslash}p{.64\linewidth}@{}}
\hline
\textbf{Choice} & \textbf{Purpose}\\\hline
\(k=\lfloor n^\kappa\rfloor,\ K=n^\tau\) & Keep \(k\) singular directions in the Gaussian replacement and frequencies up to \(K\) in the curvature calculation. The integrated replacement error tends to zero since \(3\tau<\kappa/2\).\\
\(m=\lfloor n^{4/5}\rfloor\) & Retain \(m\) tangential coordinates. The tail quadratic form has variance \(O(m^{-3})\), giving fluctuations small compared with \(1/n\).\\
\(r=n^{-1.075}\) & Smooth the retained coordinates. The noise in the width is small, while the resulting Fourier cutoff satisfies \(\sqrt m/(r\sqrt n)\le n^{0.975}\), as required by the lattice estimate.\\
\(\Delta=n^{-1.05},\quad \epsilon_n=n^{-1.02}\) & \(\Delta\) bounds the absolute change in the width. Dividing by a lower bound for that width gives a relative crossing error smaller than \(\epsilon_n\). Both errors leave room for the \(1/(2n)\) correction.\\\hline
\end{tabular}
\endgroup

\medskip\noindent\textbf{External inputs.}
The arithmetic and rectangular singular-value bounds in Proposition~\ref{pr:LCD} and Lemma~\ref{le:RV} come from Rudelson and Vershynin. We use classical concentration and convex-geometric results where cited, and Tao and Vu's theorem for the last range of thresholds. Lemma~\ref{le:rig} gives our row-minor proof of the singular-value bound needed for delocalisation. The estimates specific to the sharp coefficient are the volume and angular comparisons (Lemmas~\ref{le:zonoid}--\ref{le:angle}), the flatness estimate (Proposition~\ref{pr:flat}), and the weighted crossing argument (Proposition~\ref{pr:poly}).

\clearpage
\subsection{Notation}
We use \(|x|\) for Euclidean norm and \(\|M\|,\|M\|_F\) for the operator and Hilbert--Schmidt norms. For an \(N\times n'\) matrix with \(N\le n'\), the singular values \(s_1\le\cdots\le s_N\) are the square roots of the eigenvalues of \(MM^\top\), including zeros. The random sign row is uniform on \(\{-1,1\}^n\). Singular bases can be chosen measurably given the block. The letters \(c,C\) denote positive absolute constants whose values may change. Constants with descriptive subscripts keep their stated values.

\smallskip
\begingroup\small\renewcommand{\arraystretch}{1.17}
\noindent\begin{tabular}{@{}p{.29\linewidth}@{\hspace{.035\linewidth}}>{\raggedright\arraybackslash}p{.66\linewidth}@{}}
\hline
\textbf{Notation} & \textbf{Meaning and first use}\\\hline
\(d_j,\ a_j=\sqrt n\,d_j\) & Singular values of the block \(B\) and their scaled versions. The index runs from 1 to \(n-1\). The scalar \(a_n\) in Proposition~\ref{pr:flat} denotes the common slope.\\
\(\ell,\ y_j,\ z_j=y_j/a_j\) & Normal, tangential and rescaled tangential coordinates of the last row in the singular basis of \(B\).\\
\(W_s,\ V_s\) & \(W_s\) gives the half-width \(sW_s\) in \eqref{eq:cross}. At a crossing \(s=X_n<a_1\), the weight is \(V_s^{-2}=nu_n^2\). Stars mark the smoothed quantities.\\
\(\Acal,\ \Gcal,\ \Hcal,\ \Hcal_t\) & Successive restrictions on the block: arithmetic regularity, the five conditions in Definition~\ref{def:good}, then lower cutoffs on \(a_1\) for the two scale ranges.\\
\(\operatorname{LCD}(x),\ H\) & Least common denominator of a direction in Section~2.3, and matrix of the omitted tangential directions in \eqref{eq:H}.\\
\(D_B(s),\ \mathcal K_B,\ \bar D_n\) & Smoothed weighted crossing density, the signed integral governing its quadratic upper bound, and the averaged density at zero over retained blocks (Section~\ref{sec:poly}).\\
\(p_B^{\rm flat},\ p_B^{\rm pol}\) & Joint densities of the smoothed normal and tangential coordinates given \(B\). The normal noise has scale \(e^{-dn}\) in Section~\ref{sec:flat} and \(n^{-90}\) in Section~\ref{sec:poly}.\\
\(\omega,\ dS,\ \phi_{\sigma,d}\) & Uniform probability and surface area on \(S^{n-1}\), with \(dS=|S^{n-1}|\,d\omega\) and \(|S^{n-1}|=2\pi^{n/2}/\Gamma(n/2)\). The density of \(N(0,\sigma^2 I_d)\) is \(\phi_{\sigma,d}\), and \(\ga=(2\pi)^{-1/2}\).\\\hline
\end{tabular}
\endgroup
\par\medskip

\section{Preliminaries}
\subsection{The last row}
Fix the first \(n-1\) rows and imagine revealing the last row. If this last row lies in the span of the others, the square matrix is singular. For a nearly singular matrix, the relevant quantities are its distance from that span and its coordinates within it. The next calculation gives the exact relation between them. In this subsection let \(n\ge2\).
Let \(B=\sum_{j=1}^{n-1}d_j w_j u_j^\top\) be an \((n-1)\times n\) matrix with \(0<d_1\le\cdots\le d_{n-1}\). Here \(w_j\) and \(u_j\) are left and right singular vectors, respectively. Let \(v\) be a unit vector in the kernel of \(B\). Append \(X^\top\) to \(B\) to form \(A\), and set \(\ell=\langle v,X\rangle\), \(y_j=\langle u_j,X\rangle\). For \(0\le\lambda<d_1^2\), set
\[
 b(\lambda)=\lambda\Bigl(1+\sum_j\frac{y_j^2}{d_j^2-\lambda}\Bigr),
 \qquad b'(\lambda)=1+\sum_j\frac{d_j^2y_j^2}{(d_j^2-\lambda)^2}.
\]
\begin{lemma}\label{le:schur}
The following assertions hold.
\begin{enumerate}[label=\textup{(\roman*)}]
\item \(\smin(A)\le d_1\le s_2(A)\).
\item For \(0\le h<d_1\), \(\smin(A)\le h\) if and only if \(\ell^2\le b(h^2)\).
\item If \(\lambda=\smin(A)^2<d_1^2\), the eigenvalue \(\lambda\) of \(AA^\top\) is simple. For a unit eigenvector \(u\), its last coordinate \(u_n=\langle e_n,u\rangle\) satisfies \(u_n^2=1/b'(\lambda)\).
\item Every unit eigenvector \(u\) for the least eigenvalue satisfies
\(\ |\ell u_n|\le \smin(A)\) and \(u_n^2(1+\sum_j y_j^2/d_j^2)\le1\).
\end{enumerate}
\end{lemma}
\begin{proof}
Rank-one interlacing applied to \(A^\top A=B^\top B+XX^\top\) gives (i). In the basis \(u_1,\ldots,u_{n-1},v\), the determinant lemma gives
\[
 \det(A^\top A-\lambda I)=\prod_j(d_j^2-\lambda)\,[\ell^2-b(\lambda)],
 \qquad 0\le\lambda<d_1^2,
\]
where the identity at zero follows by continuity. Throughout this interval the product is positive. Since \(b\) increases strictly from zero, the determinant has a root in \([0,h^2]\) precisely when \(\ell^2\le b(h^2)\). This proves (ii).
Write \(u=(u',u_n)\) and \(\lambda=\smin(A)^2\). The first \(n-1\) coordinates of \(AA^\top u=\lambda u\) give
\begin{equation}\label{eq:left}
 (d_j^2-\lambda)\langle w_j,u'\rangle=-u_n d_j y_j .
\end{equation}
For \(\lambda<d_1^2\), simplicity follows from (i), and normalization in \eqref{eq:left} gives (iii). In general, \(|\langle w_j,u'\rangle|\ge |u_n y_j|/d_j\), including when \(d_j^2=\lambda\), in which case \(u_n y_j=0\). Sum the squares. Finally \(Av=\ell e_n\), so \(|\ell u_n|=|\langle A^\top u,v\rangle|\le\smin(A)\).
\end{proof}
Put \(a_j=\sqrt n\,d_j,\ z_j=y_j/a_j\), and for \(0\le s<a_1\) set
\begin{equation}\label{eq:WV}
 W_s^2=\frac1n+\sum_j\frac{z_j^2}{1-s^2/a_j^2},\qquad
 V_s^2=\frac1n+\sum_j\frac{z_j^2}{(1-s^2/a_j^2)^2}.
\end{equation}
Then \(b(s^2/n)=s^2W_s^2,\ b'(s^2/n)=nV_s^2\), and
\begin{align}
 (sW_s)'&=V_s^2/W_s\ge W_s, \label{eq:jac}\\
 X_n\le t&\ \Longleftrightarrow\ |\ell|\le tW_t \quad(0\le t<a_1),\label{eq:cross}\\
 nu_n^2&=V_{X_n}^{-2}\quad(X_n<a_1),\label{eq:weight}\\
 nu_n^2&\le \min(n,W_0^{-2}),\qquad |\ell|\sqrt n\,|u_n|\le X_n.\label{eq:crudeweight}
\end{align}
The function \(s\mapsto sW_s\) therefore increases. Its intersection with the level \(|\ell|\), when it occurs below \(a_1\), is precisely \(X_n\). At that point \eqref{eq:weight} identifies the squared last coordinate of the least left singular vector. The following identity allows us to use that coordinate as a weight.
\begin{lemma}[Exchangeability]\label{le:exchange}
Conditionally on \(A\), choose \(u\) uniformly on the unit sphere of the least eigenspace of \(AA^\top\). Then \(F_n(t)=\E[nu_n^2\one{X_n\le t}]\) for \(t\ge0\).
\end{lemma}
\begin{proof}
Permuting the rows sends this conditional uniform measure to the corresponding uniform measure for the permuted matrix. Hence \(\E[\langle e_i,u\rangle^2\one{X_n\le t}]\) is independent of \(i\). Sum over \(i\) and use \(|u|=1\).
\end{proof}

\subsection{Concentration and a Fourier majorant}
We will apply the first bound below to the length of a projected sign row and the second lemma to the probability that its normal coordinate lies in a short interval.
\begin{lemma}[Hanson--Wright]\label{le:HW}
For a sign vector \(X\), a real symmetric matrix \(M\), and \(Q_M=X^\top M X-\tr M\), one has \(\E Q_M=0\) and
\begin{align*}
 \PP(|Q_M|>s)&\le2\exp[-c_{\rm hw}\min(s^2/\|M\|_F^2,s/\|M\|)],\qquad s\ge0,\\
 \E e^{\lambda Q_M}&\le\exp(C_{\rm hw}\lambda^2\|M\|_F^2),\qquad |\lambda|\|M\|\le c_{\rm hw},\\
 \|Q_M\|_{L^p}&\le Cp\|M\|_F,\qquad p\ge1.
\end{align*}
\end{lemma}
\begin{proof}
The diagonal part cancels because \(X_i^2=1\). Removing the diagonal increases the operator norm by at most a factor of two and does not increase the Hilbert--Schmidt norm. The first two assertions follow from the tail and moment-generating function bounds in \cite{hw,versh}. Integration of the tail gives the last inequality.
\end{proof}
\begin{lemma}[Fourier majorant]\label{le:Fejer}
If \(Y\ge0\) is integrable and \(\xi,\xi'\) are real random variables, then for \(\varepsilon,\beta>0\),
\begin{align}
 \E[Y\one{|\xi|\le \varepsilon}]&\le C_\psi \varepsilon\int_{-1/\varepsilon}^{1/\varepsilon}|\E[Y e^{ik\xi}]|\,dk,\label{eq:major1}\\
 \PP(|\xi|\le \varepsilon,\ |\xi'|\le \beta)&\le C_\psi^2\varepsilon\beta\int_{-1/\varepsilon}^{1/\varepsilon}\int_{-1/\beta}^{1/\beta}
        |\E e^{i(k\xi+h\xi')}|\,dh\,dk .\label{eq:major2}
\end{align}
Here \(C_\psi=(2\sin(1/2))^{-2}\).
\end{lemma}
\begin{proof}
The function \(\psi(x)=C_\psi(\sin(x/2)/(x/2))^2
=C_\psi\int_{-1}^1(1-|k|)e^{ikx}\,dk\) is nonnegative and at least one on \([-1,1]\). Apply it to \(\xi/\varepsilon\), or take the product of two such majorants.
\end{proof}
For iid centered \(N_j\) with all moments and fixed \(p\ge2\), the Marcinkiewicz--Zygmund inequality gives
\begin{equation}\label{eq:mom}
 \Big\|\sum_j\alpha_jN_j\Big\|_{L^p}\le C_p|\alpha|,\qquad
 \Big\|\sum_j\alpha_j(N_j^2-\E N_j^2)\Big\|_{L^p}\le C_p|\alpha| .
\end{equation}
The constants depend only on \(p\) and the common law. We also use \(\|\langle x,X\rangle\|_{L^p}\le C_p|x|\).

\subsection{Arithmetic structure}
For a sign row, the characteristic function is a product of cosines. It can be large when the vector of frequencies comes close to \(\pi\Z^n\). We need to exclude such frequencies in the direction normal to \(B\), and later in nearby directions as well. Throughout the probabilistic statements below the last sign row \(X^\top\) is independent of \(B\).
Given \(\delta,\rho\in(0,1)\), a unit vector is compressible if its distance to some vector with at most \(\delta n\) nonzero coordinates is at most \(\rho\). All other unit vectors are incompressible. For \(\gamma_0,\mu\in(0,1)\), define
\[
 \operatorname{LCD}(x)=\inf\{q>0:\dist(qx,\Z^n)<\min(\gamma_0 q,\mu\sqrt n)\}.
\]
\begin{proposition}[Rudelson--Vershynin]\label{pr:LCD}
There are absolute \(\delta,\rho,\gamma_0,\mu,c_0\in(0,1)\) such that for every sufficiently large \(n\), with probability at least \(1-e^{-c_0n}\), an \((n-1)\times n\) sign matrix \(B\) has full row rank, \(|Bx|\ge c_0\sqrt n\) for all compressible \(x\), and
\(|Bx|\ge c_0 n/\operatorname{LCD}(x)\) for every incompressible \(x\) with \(\operatorname{LCD}(x)\le e^{c_0 n}\). Denote this event by \(\Acal\).
\end{proposition}
\begin{proof}
The rectangular least singular value bound of Rudelson and Vershynin \cite[Theorem 1.1]{rvrect} applied to \(B^\top\) gives full row rank except with probability \(e^{-cn}\). Their Lemma 2.6 controls compressible vectors for matrices with at least half as many rows as columns. Since \(B\) has \(n-1\ge n/2\) rows, it gives \(|Bx|\ge c\sqrt n\) simultaneously for all such \(x\), again except with probability \(e^{-cn}\).

The numbering here refers to the published version of \cite{rvrect}. For the incompressible vectors, their Section 4 treats a matrix with \(N-m\) rows and \(N\) columns. Here \(N=n\), \(m=1\), and their lattice parameters are \(\alpha=\mu\sqrt n\) and \(\gamma=\gamma_0\). Their Lemma 3.8 supplies a lower cutoff \(\operatorname{LCD}(x)\ge b_0\sqrt n\). For \(b_0\sqrt n\le\Lambda\le b_1\sqrt n\,e^{b_1 n}\), their Lemma 4.10 gives
\[
 \PP\!\left(\inf_{\substack{x\ {\rm incompressible}\\ \Lambda\le \operatorname{LCD}(x)<2\Lambda}}
              |Bx|< b_2n/\Lambda\right)\le 2e^{-n}
\]
with fixed positive \(b_0,b_1,b_2\) and a suitable \(\mu\).
Take \(\Lambda_j=2^j b_0\sqrt n\) for all \(j\ge0\) with \(\Lambda_j\le e^{b_1n/2}\). These intervals cover every admissible \(\operatorname{LCD}(x)\) from \(b_0\sqrt n\) through \(e^{b_1n/2}\). There are at most \(1+\log_2(e^{b_1n/2}/(b_0\sqrt n))=O(n)\) of them for large \(n\). Their combined failure probability is at most \(Cn e^{-n}\). On a successful interval \(|Bx|\ge b_2n/\Lambda_j\ge b_2n/\operatorname{LCD}(x)\). Choosing \(c_0\) smaller than \(b_1/2\), \(b_2\), and the absolute constants in the first paragraph proves all three assertions on the same event.
\end{proof}
For a sign vector \(X\), let \(\varphi(w)=\E e^{i\langle w,X\rangle}=\prod_i\cos w_i\).
\begin{lemma}[Lattice distance]\label{le:lattice}
Suppose \(B\in\Acal\), \(w\ne0\), and
\begin{equation}\label{eq:latticecond}
 |w|\le\pi e^{c_0n},\qquad |Bw|<c_0\sqrt n\,|w|,\qquad |Bw|<\pi c_0n .
\end{equation}
Then \(\dist(w,\pi\Z^n)\ge\min(\gamma_0|w|,\pi\mu\sqrt n)\). With \(c_1=\frac12\min(\gamma_0^2,\pi^2\mu^2)\), for every \(S\subset[n]\),
\begin{equation}\label{eq:walsh}
 \Big|\E\Big[e^{i\langle w,X\rangle}\prod_{\iota\in S}X_\iota\Big]\Big|
 \le e^{2|S|}\exp[-c_1\min(|w|^2,n)].
\end{equation}
\end{lemma}
\begin{proof}
The middle inequality makes \(x=w/|w|\) incompressible. Violation of the asserted distance bound would give \(\operatorname{LCD}(x)\le |w|/\pi\le e^{c_0n}\). Proposition~\ref{pr:LCD} would then give \(|Bw|\ge\pi c_0n\). For the second assertion use \(|\cos z|\le e^{-\dist(z,\pi\Z)^2/2}\). A factor \(X_\iota\) replaces one cosine by \(i\sin w_\iota\), losing at most \(\pi^2/8<2\) in the exponent.
\end{proof}
\begin{corollary}\label{co:slab}
For \(B\in\Acal\), the unit normal \(v\) is incompressible, \(\operatorname{LCD}(v)>e^{c_0n}\), and
\(|\varphi(kv)|\le e^{-c_1\min(k^2,n)}\) for \(|k|\le\pi e^{c_0n}\). With \(c_2=\min(c_0,c_1)\),
\begin{enumerate}[label=\textup{(\alph*)}]
\item \(\sup_z\PP(|\ell-z|\le\varepsilon\mid B)\le C(\varepsilon+e^{-c_2n})\) for \(\varepsilon\ge0\).
\item If \(\varepsilon,\beta\in[n^{-100},1]\) and \(\beta\ge a_1\), then \(\PP(|\ell|\le\varepsilon,\ |y_1|\le\beta\mid B)\le C(\varepsilon\beta+e^{-c_2n})\).
\end{enumerate}
\end{corollary}
\begin{proof}
The equation \(Bv=0\) excludes both a compressible normal and an incompressible normal with \(\operatorname{LCD}(v)\le e^{c_0n}\), by Proposition~\ref{pr:LCD}. Lemma~\ref{le:lattice} with \(w=kv\) gives the characteristic function bound. For (a), take \(\varepsilon\ge e^{-c_0n}\) in \eqref{eq:major1}. The phase from the shift by \(z\) has modulus one. Up to frequency \(\sqrt n\), the integral is bounded by the integral of \(e^{-c_1k^2}\) over the line. At larger frequencies the integrand is at most \(e^{-c_1n}\) and the integration interval has length at most \(2/\varepsilon\). This gives \(C(\varepsilon+e^{-c_1n})\). For smaller \(\varepsilon\), enlarge the slab to width \(e^{-c_0n}\).

For (b), use the frequency vector \(w=kv+hu_1\), where \(|k|\le1/\varepsilon\) and \(|h|\le1/\beta\). Orthogonality gives \(|w|^2=k^2+h^2\le2n^{200}\). Since \(a_1\le\beta\le1\),
\[
 |Bw|=|h|a_1/\sqrt n\le\min(1,|w|)/\sqrt n.
\]
Thus every nonzero \(w\) in the integration rectangle satisfies \eqref{eq:latticecond}. Inside the disk \(k^2+h^2\le n\), integrate the bound \(e^{-c_1(k^2+h^2)}\) over the plane. On the remaining rectangle use \(e^{-c_1n}\). Formula \eqref{eq:major2} gives \(C\varepsilon\beta+Ce^{-c_1n}\).
\end{proof}

\subsection{Singular values of sign matrices}
Here \(M\) is an \(N\times n'\) sign matrix with \(N\le n'\) and \(N,n'\in\{n-1,n\}\). Write \(a_j(M)=\sqrt n\,s_j(M)\) and \(a_j=a_j(B)\).
\begin{lemma}\label{le:RV}
For absolute \(c,c',C>0\) and large \(n\),
\begin{enumerate}[label=\textup{(\alph*)}]
\item \(\PP(a_1\le x)\le Cx^2+e^{-cn}\) for \(x\ge0\).
\item \(\PP(a_j(M)\le c'j)\le e^{-j}+e^{-cn}\) for \(1\le j\le N\).
\item \(\PP(\rank M<N)\le e^{-cn}\).
\end{enumerate}
\end{lemma}
\begin{proof}
We use the following form of the rectangular theorem \cite[Theorem 1.1]{rvrect}. Let \(Y\) be an \(h\times d\) matrix, \(h\ge d\), with independent copies of a centered variance-one subgaussian variable as entries. For \(\varepsilon>0\),
\[
 \PP\!\left(\smin(Y)\le\varepsilon(\sqrt h-\sqrt{d-1})\right)
     \le(C_{\rm RV}\varepsilon)^{h-d+1}+e^{-c_{\rm RV}h}.
\]
The constants depend only on the subgaussian norm of the entry distribution, hence are absolute for signs.

\step{(a) The block \(B\)} For \(Y=B^\top\) the dimensions are \(h=n\) and \(d=n-1\). The exponent is two and
\(\sqrt n-\sqrt{n-2}=2/(\sqrt n+\sqrt{n-2})\ge n^{-1/2}\).
Take \(\varepsilon=x\) for \(x>0\). The event \(a_1\le x\) is contained in the event on the left of the rectangular bound. A decreasing limit gives the same estimate at zero.

\step{(b) The \(j\)-th singular value} Delete \(j-1\) rows of \(M\) and call the remaining matrix \(M'\). Its row Gram matrix is a principal submatrix of \(MM^\top\), so interlacing gives \(s_1(M')\le s_j(M)\). Apply the rectangular theorem to \(Y=(M')^\top\) with \(h=n'\), \(d=N-j+1\). Then \(h-d+1=n'-N+j\ge j\) and
\[
 \sqrt h-\sqrt{d-1}=\frac{n'-N+j}{\sqrt{n'}+\sqrt{N-j}}\ge\frac{j}{2\sqrt n}.
\]
Choose \(\varepsilon=e^{-1}/C_{\rm RV}\) and \(c'=\varepsilon/2\). If \(a_j(M)\le c'j\), then \(s_1(M')\le\varepsilon(\sqrt h-\sqrt{d-1})\). The polynomial term in the theorem is \(e^{-(h-d+1)}\le e^{-j}\). Since \(n'\ge n-1\), the remaining term is at most \(e^{-cn}\) for an absolute \(c>0\).

\step{(c) The rank} For \(Y=M^\top\), one has \(h=n'\) and \(d=N\). Rank deficiency is precisely \(\smin(Y)=0\). Letting \(\varepsilon\downarrow0\) leaves only the exponential term.
\end{proof}
\begin{lemma}\label{le:many}
For every fixed \(\varepsilon\in(0,1)\) and \(n\ge n_0(\varepsilon)\),
\(\PP(s_{\lceil\varepsilon n/6\rceil}(M)>\varepsilon\sqrt n)\le e^{-\varepsilon^2n/80}\).
\end{lemma}
\begin{proof}
Let \(N_\varepsilon\) count the singular values at most \(\varepsilon\sqrt n\). By the Marchenko--Pastur theorem \cite[Theorem~3.6, p.~47]{bai}, \(N_\varepsilon/n\) converges to
\(\pi^{-1}\int_0^\varepsilon\sqrt{4-u^2}\,du\ge\sqrt3\,\varepsilon/\pi\).
Since \(0\le N_\varepsilon/n\le1\), this convergence gives \(\E N_\varepsilon\ge\varepsilon n/3\) for large \(n\). If \(s_{\lceil\varepsilon n/6\rceil}>\varepsilon\sqrt n\), the count is smaller than \(\varepsilon n/6\). Replacing one row changes the count by at most two, by interlacing. The bounded differences inequality \cite{mcdiarmid} thus gives an upper bound \(\exp[-2(\varepsilon n/6)^2/(4N)]\le e^{-\varepsilon^2n/72}\), which implies the claim.
\end{proof}
\begin{lemma}[Crude rigidity]\label{le:rig}
For each prescribed failure exponent \(b_{\rm fail}\ge1\), with probability at least \(1-n^{-b_{\rm fail}}\) for \(n\ge n_0(b_{\rm fail})\), \(M\) has full row rank and \(a_p(M)\le C_{\rm rig}(b_{\rm fail})p^2\log n\) for \(1\le p\le n^{1/4}\).
\end{lemma}
\begin{proof}
To bound \(a_p(M)\) from above, we will bound the inverse-square sum \(\sum_{j\ge p}a_j(M)^{-2}\) from below. Lemma~\ref{le:RV}(b) controls the contribution far out in this sum. The remaining indices must then carry a definite amount of mass. We obtain that mass by measuring distances of individual rows from spans of other rows.

\step{A deterministic comparison of minors} Write \(R_i\) for the rows and \(\mathfrak g(I)\) for the Gram determinant of \((R_i)_{i\in I}\), with \(\mathfrak g(\varnothing)=1\). On full row rank, put \(x_j=s_j(M)^{-2}\) in decreasing order. Set \(e_0=1\). The elementary symmetric polynomial \(e_p=e_p(x)\) is the sum of the principal \(p\times p\) minors of \((MM^\top)^{-1}\). Jacobi's complementary-minor formula gives
\begin{equation}\label{eq:minor}
 e_p=\sum_{|S|=p}\frac{\mathfrak g(S^c)}{\mathfrak g([N])}.
\end{equation}
For \(i\in S\), let \(\operatorname{dist}_i(S)\) be the distance of \(R_i\) to the span of the rows indexed by \(S^c\). Since \(\mathfrak g(S^c\cup\{i\})=\operatorname{dist}_i(S)^2\mathfrak g(S^c)\),
\begin{equation}\label{eq:minorlower}
 p e_p=\sum_{|S|=p}\sum_{i\in S}\frac{\mathfrak g((S\setminus\{i\})^c)}
                  {\operatorname{dist}_i(S)^2\mathfrak g([N])}
 \ \ge\ \frac{N-p+1}{\max_{S,i}\operatorname{dist}_i(S)^2}\,e_{p-1}.
\end{equation}
Indeed, each \((p-1)\)-set \(S\setminus\{i\}\) arises \(N-p+1\) times. For an upper bound, group the monomials in \(e_p\) by their largest index \(j\). Such a monomial contributes \(x_j\) times a monomial of \(e_{p-1}\), and \(j\ge p\). Hence
\begin{equation}\label{eq:minorupper}
 e_p\le e_{p-1}\sum_{j\ge p}x_j.
\end{equation}
\step{Distances of independent rows} For fixed \(S\) and \(i\in S\), condition only on the rows in \(S^c\). Whenever these rows are linearly independent, the projection onto their orthogonal complement has rank \(p+n'-N\le2p\). The remaining row \(R_i\) is independent of this projection. Lemma~\ref{le:HW} bounds the conditional probability of \(\operatorname{dist}_i(S)^2>C_3p\log n\) by \(2e^{-c C_3p\log n}\). Full row rank implies linear independence of every such subset. A union bound over at most \(n^{p+1}\) pairs at level \(p\), summed over \(p\le n^{1/4}\), costs at most \(n^{-b_{\rm fail}-1}\) for \(C_3=C_3(b_{\rm fail})\) large. Lemma~\ref{le:RV} also gives full rank and \(a_j(M)\ge c'j\) for \(j\ge J=\lceil(b_{\rm fail}+3)\log n\rceil\), except with probability at most \(n^{-b_{\rm fail}-1}\) for large \(n\).
\step{Locating the inverse-square mass} On the event just constructed, divide \eqref{eq:minorlower} by \(p e_{p-1}\) and compare with \eqref{eq:minorupper}. Since \(a_j(M)^{-2}=x_j/n\) and \(N-p+1\ge n/2\) for large \(n\), the result is
\(\sum_{j\ge p}a_j(M)^{-2}\ge(2C_3p^2\log n)^{-1}\).
With \(J'=J+\lceil8 C_3(c')^{-2}p^2\log n\rceil\), the contribution from \(j\ge J'\) is at most \(2(c')^{-2}/J'\le(4C_3p^2\log n)^{-1}\).
Subtracting leaves
\[
 \frac{1}{4C_3p^2\log n}\le \sum_{p\le j<J'}a_j(M)^{-2}\le J'a_p(M)^{-2}.
\]
Here \(J'\le C_{\rm rig}(b_{\rm fail}) p^2\log n\), so taking square roots gives \(a_p(M)\le C_{\rm rig}(b_{\rm fail}) p^2\log n\) after enlarging \(C_{\rm rig}(b_{\rm fail})\).
\end{proof}
\begin{lemma}[Delocalisation]\label{le:deloc}
For each prescribed failure exponent \(b_{\rm fail}\ge1\) and \(n\ge n_0(b_{\rm fail})\), with probability at least \(1-n^{-b_{\rm fail}}\) every unit \(w\) satisfying \(|Bw|\le n^{-1/4}\) has \(\|w\|_\infty\le n^{-1/4}\).
\end{lemma}
\begin{proof}
Delete column \(b_i\) to get \(B^{(i)}\). Let \(\Pi_i\) project onto the \(p=\lfloor n^{1/8}\rfloor\) least left singular directions of this square minor. By Lemma~\ref{le:rig}, \(\|\Pi_iB^{(i)}\|\le C_{\rm rig}(b_{\rm fail}+2)n^{-1/4}\log n\) except with probability \(n^{-b_{\rm fail}-2}\). By independence and Lemma~\ref{le:HW}, \(|\Pi_i b_i|\ge\sqrt{p/2}\) except with probability \(2e^{-cp}\). Project the decomposition \(Bw=b_iw_i+B^{(i)}w_{-i}\) onto this same subspace, where \(w_{-i}\) denotes \(w\) with coordinate \(i\) removed. Since \(|w_{-i}|\le1\),
\[
 |w_i|\sqrt{p/2}\le |Bw|+\|\Pi_iB^{(i)}\|.
\]
Take the union over \(i\) and use \(n^{-1/16}\log n\to0\).
\end{proof}
\begin{lemma}[Approximate null vectors]\label{le:flatnull}
For every fixed \(a\in(0,1)\) there is \(\beta_a>0\) such that, for sufficiently large \(n\), with probability at least \(1-e^{-\beta_a n}\) every unit \(v\) satisfying \(|Av|\le\beta_a\sqrt n\) has \(\|v\|_\infty\le a\).
\end{lemma}
\begin{proof}
Take \(\varepsilon=a^2/48\) and \(p=\lceil\varepsilon n/6\rceil\), so \(2\le p\le n-1\) for large \(n\). Delete column \(A_j\) from \(A\), giving an \(n\times(n-1)\) matrix \(A^{(j)}\). Let \(\Pi_j\) project onto the bottom \(p\) eigendirections of \(A^{(j)}(A^{(j)})^\top\). Its spectrum consists of the squared singular values of \((A^{(j)})^\top\) and one additional zero, counting multiplicities. Consequently
\[
 \|\Pi_j A^{(j)}\|=s_{p-1}((A^{(j)})^\top)\le s_p((A^{(j)})^\top).
\]
Lemma~\ref{le:many}, applied to the \((n-1)\times n\) transpose, bounds the right side by \(\varepsilon\sqrt n\) except with probability \(e^{-\varepsilon^2n/80}\).

Given the minor, \(A_j\) is an independent sign column and \(\Pi_j\) is a fixed rank-\(p\) projection. Lemma~\ref{le:HW} yields \(|\Pi_j A_j|\ge\sqrt{p/2}\) except with conditional probability \(2e^{-cp}\). Projecting the identity \(Av=A_jv_j+A^{(j)}v_{-j}\) now gives
\[
 |v_j|\sqrt{p/2}\le |Av|+\varepsilon\sqrt n .
\]
Thus \(|Av|\le\varepsilon\sqrt n\) implies \(|v_j|\le2\sqrt{12\varepsilon}=a\). Take a union bound over the \(n\) columns and choose \(\beta_a\le\varepsilon\) small enough to absorb both failure probabilities.
\end{proof}

\subsection{Parameters and a good event}
Fix \(0<\theta\le1/100\) as specified in Section~\ref{sec:finish}. Set
\begin{equation}\label{eq:params}
\begin{gathered}
 \kappa=\theta/20,\quad k=\lfloor n^\kappa\rfloor,\quad m=\lfloor n^{4/5}\rfloor,\quad
 r=n^{-1.075},\quad M_0=n^6,\\
 \tau=\theta/1000,\qquad K=n^\tau,\qquad J_0=\lceil103\log n\rceil,\qquad \Delta=n^{-1.05},\quad\epsilon_n=n^{-1.02}.
\end{gathered}
\end{equation}
\begin{definition}\label{def:good}
Write \(C_*=C_{\rm rig}(102)\) from Lemma~\ref{le:rig}. The event \(B\in\Gcal\) consists of the following five conditions.
\begin{enumerate}[label=(G\arabic*)]
\item \(B\in\Acal\).
\item \(a_j\ge c'j\) for \(J_0\le j\le n-1\).
\item \(a_j\le C_* j^2\log n\) for \(1\le j\le n^{1/4}\).
\item \(a_m\le c_0 n/2\).
\item Every unit \(w\in\operatorname{span}(v,u_1,\ldots,u_k)\) satisfies \(\|w\|_\infty\le n^{-1/4}\).
\end{enumerate}
\end{definition}
\begin{proposition}\label{pr:good}
For large \(n\), \(\PP(B\notin\Gcal)\le n^{-100}\).
\end{proposition}
\begin{proof}
Use Proposition~\ref{pr:LCD} for (G1), sum Lemma~\ref{le:RV}(b) over \(j\ge J_0\) for (G2), and apply Lemma~\ref{le:rig} for (G3). For (G4), use Lemma~\ref{le:many} with \(\varepsilon=c_0/2\), since \(m\le\lceil c_0n/12\rceil\) eventually. On (G3), unit vectors in the stated span satisfy \(|Bw|\le a_k/\sqrt n\le C_* k^2 n^{-1/2}\log n\le n^{-1/4}\), so Lemma~\ref{le:deloc} yields (G5).
\end{proof}
On \(B\in\Gcal\) we will keep the first \(m\) tangential coordinates and estimate the contribution of the rest through the matrix
\begin{equation}\label{eq:H}
 H=\sum_{j>m}a_j^{-2}u_ju_j^\top.
\end{equation}
For large \(n\), \(m\ge J_0\), so (G2) applies to every eigenvalue \(a_j^{-2}\) in this sum. For any integer \(q\ge J_0\), it also applies to every \(j>q\). Comparing \(\sum_{j>q}j^{-2}\) and \(\sum_{j>m}j^{-4}\) with their integrals gives the first four bounds below. For the last bound split at \(J_0\) and use \(a_j\ge a_1\) on the first \(J_0-1\) terms.
\begin{equation}\label{eq:tails}
\begin{gathered}
 \|H\|\le Cm^{-2},\quad \tr H\le Cm^{-1},\quad \tr H^2\le Cm^{-3},\\
 \sum_{j>q}a_j^{-2}\le C/q,\qquad \sum_{j\ge1}a_j^{-2}\le C(1+a_1^{-2}\log n).
\end{gathered}
\end{equation}
Lemma~\ref{le:HW} applied to the projection onto the first \(k\) directions gives \(\PP(E_k\mid B)\le2e^{-ck}\), where
\begin{equation}\label{eq:Ek}
 E_k=\Big\{\sum_{j\le k}y_j^2<k/2\Big\}.
\end{equation}
On \(E_k^c\),
\begin{equation}\label{eq:w0}
 W_0^2\ge\sum_{j\le k}z_j^2\ge k/(2a_k^2)\ge w_0^2,\qquad
 w_0=(2C_* k^{3/2}\log n)^{-1}\ge n^{-2\kappa}
\end{equation}
for sufficiently large \(n\).
\section{A one-sided estimate at a single scale}\label{sec:area}
Fix \(t_0=n^{-25}\) and \(\sigma=t_0/\sqrt n\). We test the distribution of \(X_n\) against a Gaussian kernel concentrated near zero. The normalization gives exactly \(a\) when the kernel is integrated against \(a\,dx\) on the positive half-line, since \(\sqrt{2/\pi}\int_0^\infty e^{-x^2/2}\,dx=1\). Section~\ref{sec:flat} will give a precise version of this observation for the discrete sign ensemble.

Two one-dimensional densities enter the comparison. For \(v\in\sphere\), put
\[
 q(v)=\E\phi_{\sigma,1}(\langle v,X\rangle),\qquad
 f_n=\frac{\Gamma(n/2)}{\sqrt{\pi n}\,\Gamma((n-1)/2)} .
\]
Thus \(q(v)\) is the density at zero of \(\langle v,X\rangle+\sigma g\) for an independent standard Gaussian \(g\). The constant \(f_n\) is the density at zero of one coordinate of a uniform point of \(\sqrt n\,\sphere\). Let \(J_A(v)\) be the tangential Jacobian of \(v\mapsto Av\) on \(\sphere\). In any orthonormal basis of \(v^\perp\), form the \(n\times(n-1)\) matrix \(Y\) of projected rows of \(A\). Then \(J_A(v)=\sqrt{\det(Y^\top Y)}\).
\begin{proposition}\label{pr:area}
For all sufficiently large \(n\),
\[
 \frac{\sqrt{2/\pi}}{t_0}\E e^{-X_n^2/(2t_0^2)}
       \le 1-\frac1{2n}+Cn^{-3/2}.
\]
\end{proposition}

The proof separates two effects. For a fixed direction \(v\), independence of the rows produces the factor \(q(v)^n\). The remaining factor is an expected tangential Jacobian. Lemma~\ref{le:zonoid} bounds this Jacobian by the corresponding one for spherical rows. After normalization, the two bounds are
\[
 \underbrace{1-\frac{5}{4n}+O(n^{-2})}_{\text{spherical model}}
 \qquad\text{and}\qquad
 \underbrace{1+\frac{3}{4n}+O(n^{-3/2})}_{\text{change of row law}} .
\]
The gain in Proposition~\ref{pr:area} is the first-order term in their product. We first show why a surface integral detects a small singular value.

\subsection{Gaussian area of a thin ellipsoid}
\begin{lemma}\label{le:area}
Let \(A\) be a real \(n\times n\) matrix with singular values \(s_1\le\cdots\le s_n\), where \(s_1\le\sigma\sqrt n\) and \(s_2\ge4\sigma\sqrt n\). For \(\mathcal C_A=\{v\in\sphere:|Av|\le3\sigma\sqrt n\}\),
\[
 \int_{\mathcal C_A}\phi_{\sigma,n}(Av)J_A(v)\,dS(v)
       \ge 2(1-e^{-n/2})\,\phi_{\sigma,1}(s_1).
\]
\end{lemma}
\begin{proof}
By orthogonal invariance take \(A=\diag(s_i)\). The image of the sphere is an ellipsoid with shortest semiaxis \(s_1\). Near the two ends of that axis its boundary can be written as two graphs over the same ball in the other \(n-1\) coordinates. The area formula \cite[Corollary~3.2.20]{federer} gives
\begin{equation}\label{eq:area}
 \int_{\mathcal C_A}\phi_{\sigma,n}(Av)J_A(v)\,dS(v)
 =\int_{\R^n}\phi_{\sigma,n}(y)\,\#\{v\in \mathcal C_A:Av=y\}\,d\mathcal H^{n-1}(y).
\end{equation}
For \(|x|\le2\sigma\sqrt n\) in the last \(n-1\) coordinates, set
\(h(x)=(1-\sum_{j\ge2}x_j^2/s_j^2)^{1/2}\ge\sqrt3/2\).
The points \((\pm h(x),x_2/s_2,\ldots,x_n/s_n)\) belong to \(\mathcal C_A\), since their images have squared norm at most \(s_1^2+|x|^2\le5\sigma^2n\). Their images form two graphs \((\pm s_1h(x),x)\) over the ball. Their surface elements dominate \(dx\), and their Gaussian densities are at least \(\phi_{\sigma,1}(s_1)\phi_{\sigma,n-1}(x)\). When \(s_1=0\) the graphs coincide and the multiplicity in \eqref{eq:area} is at least two. Integrating over the ball gives at least \(2\phi_{\sigma,1}(s_1)\PP(\chi^2_{n-1}\le4n)\). The required chi-square tail estimate (see also \cite[eq.~(4.3) after Lemma~1, p.~1325]{laurent}) follows by applying Markov's inequality at parameter \(1/4\),
\[
 \PP(\chi^2_{n-1}>4n)\le e^{-n}\E e^{\chi^2_{n-1}/4}
 =e^{-n}2^{(n-1)/2}\le e^{-n/2}.
\]
\end{proof}

\subsection{A volume comparison}
\begin{lemma}[Zonoid comparison]\label{le:zonoid}
Let \(Y\) have \(d+1\) independent rows from a common law supported on the ball of radius \(r_{\rm ball}>0\) in \(\R^d\). Let \(Y_*\) have independent rows uniform on \(r_{\rm ball} S^{d-1}\). Then
\[
 \E\sqrt{\det(Y^\top Y)}\le \E\sqrt{\det(Y_*^\top Y_*)}.
\]
\end{lemma}
\begin{proof}
\step{Turn the Jacobian into a volume} Append an independent standard Gaussian column \(g\in\R^{d+1}\). For a full-rank \(Y\), the absolute determinant of \([Y,g]\) is the \(d\)-dimensional volume of the column parallelepiped times the absolute normal component of \(g\). This component has law \(|N(0,1)|\) conditionally on \(Y\), so
\begin{equation}\label{eq:append}
 \sqrt{\det(Y^\top Y)}=\sqrt{\pi/2}\,\E_g|\det[Y,g]|.
\end{equation}
Both sides vanish when \(Y\) has deficient column rank. Let \(\xi=(Y_1,g_1)\) be one row of the augmented matrix. The centered segment \([-\xi/2,\xi/2]\) has support function \(|\langle x,\xi\rangle|/2\), the largest scalar product with \(x\) over the segment. Averaging these support functions defines a convex body \(K\), the zonoid of the row law. Vitale's determinant formula \cite{vitale} gives
\begin{equation}\label{eq:vitale}
 \E|\det[Y,g]|=(d+1)!\vol_{d+1}(K).
\end{equation}
The centered segment convention for \(h_K\) gives the same volume as the convention using segments \([0,\xi]\), since the two mean bodies differ by a translation.

\step{Average over row directions} Rotate the first \(d\) coordinates by a uniform orthogonal transformation, leaving the last coordinate fixed. The support function of the averaged body \(\bar K\) is the average of the rotated support functions. The Brunn--Minkowski inequality \cite[\S7.1]{schneider} gives \(\vol(\bar K)^{1/(d+1)}\ge\vol(K)^{1/(d+1)}\), since all rotated bodies have the same volume. Indeed, finite averages of independent Haar rotations have this property, and their support functions converge uniformly almost surely by the strong law on a countable dense set and their common Lipschitz bound. For \(x=(x',x_{d+1})\),
\[
 h_{\bar K}(x)=\tfrac12\E\big|\,|Y_1|\langle x',\Theta\rangle+x_{d+1}g_1\big|,
\]
where \(\Theta\) is uniform on \(S^{d-1}\) and independent of \((Y_1,g_1)\). \step{Increase the row length} For fixed \(z\), \(r'\mapsto\E_\Theta|r'\langle x',\Theta\rangle+z|\) is even and convex, hence nondecreasing on \([0,\infty)\). Since \(|Y_1|\le r_{\rm ball}\), this support function is bounded above by that of the body for \(Y_*\). Inclusion, \eqref{eq:append}, and \eqref{eq:vitale} give the claim.
\end{proof}

\subsection{Spherical calibration}
Here \(Y_*\) has size \(n\times(n-1)\) and independent rows uniform on \(\sqrt n\,S^{n-2}\). Put \(J_*=\sqrt{\det(Y_*^\top Y_*)}\).
\begin{lemma}\label{le:cal}
For all sufficiently large \(n\),
\[
 c_n^{\sph}:=\frac{|\sphere|}{\sqrt n}\, f_n^{\,n}\E J_*
       =1-\frac{5}{4n}+O(n^{-2}).
\]
\end{lemma}
Proposition~\ref{pr:sph} proves this normalization by computing the same quantity in two ways. Lemma~\ref{le:sphgeom} identifies \(c_n^{\sph}\) with the right derivative at zero of the distribution function of \(\sqrt n\,\smin(A_{\sph})\). Lemma~\ref{le:sphder} computes that derivative from an exact formula for its distribution.

\subsection{The angular moment}
\begin{lemma}\label{le:angle}
There are absolute \(a\in(0,1/2]\) and \(C_{\rm lat}>0\) such that, for \(\Lambda_{\rm lat}=\sigma^{-1}\sqrt{C_{\rm lat}\log n}\) and
\[
 \mathcal V=\{v\in\sphere:\|v\|_\infty\le a,\ \operatorname{LCD}(v)>\Lambda_{\rm lat}\},
\]
for all sufficiently large \(n\) one has both \(q(v)\le\ga\) on \(\mathcal V\) and
\[
 \int_{\mathcal V}(q(v)/f_n)^n\,d\omega(v)\le1+\frac{3}{4n}+Cn^{-3/2}.
\]
\end{lemma}
\begin{proof}
The power \(n\) makes a uniform \(O(1/n)\) approximation to \(q/f_n\) too crude. We use the uniform bound only to discard atypical directions. For typical directions we expand the power and integrate its linear and quadratic terms separately.

Fourier inversion gives
\begin{equation}\label{eq:qFourier}
 2\pi q(v)=\int_{\R}e^{-\sigma^2 u^2/2}\prod_i\cos(uv_i)\,du .
\end{equation}
For \(v\in \mathcal V\), the definition of \(\operatorname{LCD}\) gives
\begin{equation}\label{eq:qtail}
 \Big|\prod_i\cos(uv_i)\Big|
 \le\exp[-\tfrac12\min(\gamma_0^2 u^2,\pi^2\mu^2 n)]
 \qquad(|u|\le\pi \Lambda_{\rm lat}).
\end{equation}
Choose \(C_{\rm lat}\) large enough that, for all sufficiently large \(n\),
\[
 \int_{|u|>\pi \Lambda_{\rm lat}}e^{-\sigma^2 u^2/2}\,du\le n^{-44},
 \qquad 2\pi \Lambda_{\rm lat} e^{-\pi^2\mu^2 n/2}\le n^{-44}.
\]

\step{Uniform bound} Write \(b=\|v\|_\infty\). For even \(j\) let \(m_j(x)=\sum_i x_i^j\), and abbreviate \(m_j=m_j(v)\) when the direction is fixed. For \(|u|\le1/b\), the product is positive and at most \(\exp(-u^2/2-u^4m_4/12)\), by the expansion of \(-\log\cos\). As \(1/b\ge2\), integrating on this interval gives at most \(\sqrt{2\pi}-c m_4\), by considering \(1\le |u|\le2\). For \(1/b<|u|\le\pi \Lambda_{\rm lat}\), \eqref{eq:qtail} bounds the integral by \(C_{\gamma_0}e^{-\gamma_0^2/(2b^2)}+n^{-44}\). Fix \(a\) so small that its first term is at most \(cb^4/2\) for \(0<b\le a\). Since \(m_4\ge\max(b^4,1/n)\), \eqref{eq:qFourier} gives \(q(v)\le\ga\). Stirling's formula yields
\begin{equation}\label{eq:gamma}
 f_n/\ga=1-\frac{3}{4n}+O(n^{-2}).
\end{equation}
Thus \((q/f_n)^n\le e\) on \(\mathcal V\) for large \(n\).

\step{The angular mean} Let \(\mathcal T_n=\{\|v\|_\infty\le n^{-1/2}\log n\}\). The even moments of a spherical coordinate are bounded by those of \(N(0,1/n)\). The exponential moment bound and a union bound give
\(\omega(\mathcal T_n^c)\le 2n e^{-(\log n)^2/2}\).
Write \(q_<(v)\) for the integral in \eqref{eq:qFourier}, divided by \(2\pi\), restricted to \(|u|\le\Lambda_{\rm low}=n^{1/3}\). This cutoff lies within the range \(|u|\le\pi\Lambda_{\rm lat}\) where \eqref{eq:qtail} applies.
On \(\mathcal V\), \eqref{eq:qtail} and the Gaussian cutoff give \(|q-q_<|\le Cn^{-44}\).
On all of \(\mathcal T_n\), \(\Lambda_{\rm low}\|v\|_\infty<1\), so the product of cosines on this interval is positive and at most \(e^{-u^2/2}\). In particular \(0\le q_<\le\ga\) there. Applying the mean value theorem to the \(n\)-th power on \(\mathcal V\cap \mathcal T_n\), using \eqref{eq:gamma}, gives
\begin{equation}\label{eq:angleint}
 \int_{\mathcal V}(q/f_n)^n\,d\omega
 \le \int_{\mathcal T_n}(q_</f_n)^n\,d\omega+O(n^{-10}).
\end{equation}

The purpose of the cutoff is that \(q_<\) is uniformly controlled on typical directions, including those whose large Fourier frequencies have arithmetic peaks. We can now use the fixed length \(|X|=\sqrt n\) to control the linear term. For every fixed sign row, its projection along a uniform direction has density
\(f(x)=f_n(1-x^2/n)_+^{(n-3)/2}\). Thus
\[
 \int_{\sphere}q\,d\omega=\int_{\R}\phi_{\sigma,1}(x)f(x)\,dx\le f_n,
 \qquad \int_{\sphere}\varphi(uv)\,d\omega(v)=\widehat f(u).
\]
For large \(n\), differentiating this density thirteen times gives
\( |f^{(13)}(x)|\le C(1+|x|^{13})e^{-x^2/4}\).
Indeed, the coefficients produced by differentiating \((1-x^2/n)^{(n-3)/2}\) are bounded uniformly, and every term on its support contains at least the power \((n-3)/2-13\) of \(1-x^2/n\). For \(n\ge58\) this exponent is at least \(n/4\), so the power is bounded by \(e^{-x^2/4}\). The lower derivatives vanish at the endpoints. Integration by parts therefore gives \(|\widehat f(u)|\le C|u|^{-13}\), and
\[
 \left|\int_{\sphere}(q-q_<)\,d\omega\right|
 \le C\int_{\Lambda_{\rm low}}^\infty u^{-13}\,du=O(n^{-4}).
\]
Globally \(|q_<|\le\Lambda_{\rm low}/\pi\). The bound on \(\omega(\mathcal T_n^c)\) then shows that \(h=q_</f_n-1\) satisfies \(\int_{\mathcal T_n}h\,d\omega\le Cn^{-4}\). This estimate controls the linear term in the angular moment.

\step{The angular variance} On \(\mathcal T_n\), expand \(\log\cos x=-x^2/2-x^4/12+O(x^6)\) for \(|u|\le(\log n)^2\) in \(q_<\). Here \(m_4\le n^{-1}\log^2n\) and \(m_6\le n^{-2}\log^4 n\). After extracting \(e^{-u^2/2}\), the error in replacing the product by \(1-u^4m_4/12\) is at most \(C(u^6m_6+u^8m_4^2)\). The rest of the truncated integral is bounded by the Gaussian tail, and removing \(e^{-\sigma^2u^2/2}\) costs \(O(\sigma^2)\). The fourth Gaussian moment is three, so \(q_</\ga=1-m_4/4+O(m_6+m_4^2+n^{-4})\).

The coefficient comes from the fluctuation of the fourth power sum. For a uniform direction, writing \(d_{8,n}=n(n+2)(n+4)(n+6)\), the Gaussian representation of the sphere gives \(\E_\omega v_1^4=3/[n(n+2)]\), \(\E_\omega v_1^8=105/d_{8,n}\), and \(\E_\omega v_1^4v_2^4=9/d_{8,n}\). Summing the diagonal and off-diagonal terms gives
\[
 \E_\omega m_4=\frac3{n+2},\qquad
 \operatorname{Var}_\omega(m_4)=
 \frac{24(n-1)}{(n+2)^2(n+4)(n+6)}.
\]
For fixed \(p\ge2\), one also has
\begin{equation}\label{eq:sphmom}
 \|m_4-\E_\omega m_4\|_{L^p(\omega)}\le C_p n^{-3/2},\qquad
 \|m_6\|_{L^p(\omega)}+\|m_4^2\|_{L^p(\omega)}\le C_p n^{-2}.
\end{equation}
Indeed, represent \(v=g/|g|\) for \(g\sim N(0,I_n)\). On the event \(|g|^2\ge n/2\),
\(|m_4-3/n|\le C(|m_4(g)-3n|+|m_2(g)-n|)/n^2\) and \(m_6\le8m_6(g)/n^3\).
Apply \eqref{eq:mom} to centered Gaussian powers and use \(m_4,m_6\le1\) on the exponentially unlikely complement.

By \eqref{eq:gamma} and the quartic expansion,
\[
 h=-\tfrac14(m_4-\E_\omega m_4)+e_{\rm ang}(v),\qquad
 |e_{\rm ang}(v)|\le C(m_6+m_4^2+n^{-2})\quad\text{on }\mathcal T_n.
\]
Cauchy--Schwarz and \eqref{eq:sphmom} give
\[
 \int_{\mathcal T_n}h^2\,d\omega\le\tfrac1{16}\operatorname{Var}_\omega(m_4)+O(n^{-7/2}),
 \qquad \int_{\mathcal T_n}|h|^3\,d\omega=O(n^{-9/2}).
\]
Since \(-1\le h\le C/n\) on \(\mathcal T_n\), the remainder after the quadratic term of \((1+h)^n\) is bounded by \(Cn^3|h|^3\). The constant term integrates to at most one and the linear term to at most \(Cn^{-3}\). The quadratic term contributes at most
\(\binom n2\operatorname{Var}_\omega(m_4)/16+O(n^{-3/2})=3/(4n)+O(n^{-3/2})\).
Substitution in \eqref{eq:angleint} proves the claim.
\end{proof}

\subsection{Proof of Proposition~\ref{pr:area}}
\step{Restrict the directions of integration} Let \(\mathcal E_{\rm dir}\) be the intersection of \(B\in\Acal\) with the event in Lemma~\ref{le:flatnull} for the \(a\) of Lemma~\ref{le:angle}. Its complement has probability \(O(e^{-cn})\). On \(\mathcal E_{\rm dir}\), if \(|Av|\le3t_0\), then \(\|v\|_\infty\le a\) and \(v\) is incompressible since \(|Bv|\le3t_0\). If \(\operatorname{LCD}(v)\le \Lambda_{\rm lat}\), Proposition~\ref{pr:LCD} would give \(|Bv|\ge c_0 n/\Lambda_{\rm lat}=c_0\sqrt n\,t_0/\sqrt{C_{\rm lat}\log n}>3t_0\). Thus every such \(v\) belongs to \(\mathcal V\).

Write \(s_1,s_2\) for the two least singular values of \(A\). The quantity inside the expectation in Proposition~\ref{pr:area}, including its prefactor, is \(2\phi_{\sigma,1}(s_1)/\sqrt n\). For \(s_1>\sigma\sqrt n\) it is \(O(t_0^{-1}e^{-n/2})\). On \(\mathcal E_{\rm dir}\), \(s_2\ge d_1\), so \(\PP(s_2<4t_0,\mathcal E_{\rm dir})\le Cnt_0^2+e^{-cn}\) by Lemma~\ref{le:RV}(a). The integrand on either exceptional event is bounded by \(C/t_0\). Thus the costs of \(\mathcal E_{\rm dir}^c\) and of \(s_2<4t_0\) are together \(O(nt_0+e^{-cn}/t_0)\). On the remaining event apply Lemma~\ref{le:area}, enlarge its domain to \(\mathcal V\), and then drop the event from the nonnegative surface integral. The domain \(\mathcal V\) is deterministic, so Tonelli's theorem gives
\begin{equation}\label{eq:surface}
 \frac{\sqrt{2/\pi}}{t_0}\E e^{-X_n^2/(2t_0^2)}
 \le Cnt_0+\frac{1+2e^{-n/2}}{\sqrt n}
       \int_{\mathcal V} \E[\phi_{\sigma,n}(Av)J_A(v)]\,dS(v).
\end{equation}
\step{Factor the integrand} For fixed \(v\), the Gaussian weight in \eqref{eq:surface} is the product of the \(n\) weights \(\phi_{\sigma,1}(\langle R_i,v\rangle)\), one for each sign row. Each weight has mean \(q(v)\). After dividing by these means the rows are still independent, with the same tilted law. Write \(\E_v\) for expectation under this law and \(Y\) for the matrix of rows projected onto \(v^\perp\). Then
\[
 \E[\phi_{\sigma,n}(Av)J_A(v)]
   =q(v)^n\,\E_v\sqrt{\det(Y^\top Y)} .
\] Each projected row has squared length \(n-\langle R_i,v\rangle^2\le n\). Lemma~\ref{le:zonoid} therefore bounds the last expectation by \(\E J_*\). Since \(dS=|\sphere|\,d\omega\), the right side of \eqref{eq:surface} is at most
\[
 Cn^{-24}+(1+2e^{-n/2})c_n^{\sph}\int_{\mathcal V}(q/f_n)^n\,d\omega.
\]
Lemmas~\ref{le:cal} and \ref{le:angle} give \(1+(-5/4+3/4)/n+O(n^{-3/2})\), as claimed.

\section{Flatness at very small scales}\label{sec:flat}
Section~\ref{sec:area} controls a Gaussian average at the scale \(t_0=n^{-25}\).
To extract a bound on \(F_n(t)\) from this average, we show that the ratio \(F_n(t)/t\) varies little throughout an interval containing \(t_0\).
For a fixed block \(B\), the crossing criterion \eqref{eq:cross} asks whether the normal coordinate \(\ell\) lies in an interval of half-width \(tW_t\). Both \(\ell\) and \(W_t\) depend on the last row. We smooth the normal coordinate together with the tangential projections that determine the width. At fixed smoothed tangential coordinates the interval is now known. Integrating its joint density over the normal coordinate produces a factor of \(t\), with a coefficient independent of \(t\).

\begin{proposition}\label{pr:flat}
There are absolute \(\cflat,C>0\) and, for each sufficiently large \(n\), a number \(a_n\in[0,C]\) such that
\[
 |F_n(t)-a_nt|\le Ct n^{-1.05}\qquad(e^{-\cflat n}\le t\le n^{-20}).
\]
\end{proposition}
\begin{corollary}\label{co:slope}
For large \(n\), \(a_n\le1-\frac1{2n}+Cn^{-1.05}\).
\end{corollary}
\begin{proof}
For \(t_0=n^{-25}\), integrate the identity \(e^{-y^2/2}=\int_y^\infty x e^{-x^2/2}\,dx\) at \(y=X_n/t_0\). Tonelli's theorem gives
\[
 \E e^{-X_n^2/(2t_0^2)}=\int_0^\infty x e^{-x^2/2}F_n(t_0x)\,dx ,
\]
including any atom at zero. On \(n^{-1}\le x\le\sqrt{12\log n}\) the arguments satisfy
\[
 t_0x\in[n^{-26},n^{-25}\sqrt{12\log n}]\subset[e^{-\cflat n},n^{-20}]
\]
for large \(n\). Write \(b_n=a_n-Cn^{-1.05}\), using the error constant of Proposition~\ref{pr:flat}. If \(b_n\le0\) the claim follows at once. Otherwise \(F_n(t_0x)\ge b_n t_0 x\) on this interval. The omitted mass of \(\sqrt{2/\pi}\,x^2e^{-x^2/2}\) is \(O(n^{-3})\) by integration of \(x^2\) at the lower end and the Gaussian tail bound at the upper end. Proposition~\ref{pr:area} thus gives
\[
 b_n(1-O(n^{-3}))\le \frac{\sqrt{2/\pi}}{t_0}\E e^{-X_n^2/(2t_0^2)}
 \le 1-\frac{1}{2n}+O(n^{-3/2}).
\]
This yields the stated bound on \(a_n\).
\end{proof}

\subsection{Rare events on a thin slab}
The squared width \(W_0^2\) contains the tail quadratic form \(X^\top H X\), whose conditional mean is \(\tr H\). We shall replace it by its mean and use the lower bound \(W_0\ge w_0\) from \eqref{eq:w0}. The two exceptional events are \(E_k\) from \eqref{eq:Ek} and
\begin{equation}\label{eq:EH}
 E_H=\{|X^\top H X-\tr H|>n^{-1.1}\}.
\end{equation}
An unconditional error \(e^{-n^c}\) with \(0<c<1\) can exceed the entire crossing probability at the lower endpoint. A failed width approximation affects the answer only when the normal coordinate is already in a thin slab. The following estimate retains the slab width in the error.
\begin{lemma}\label{le:rare}
For \(B\in\Gcal\), \(e^{-c_0n}\le\varepsilon\le1\), and sufficiently large \(n\),
\[
 \PP(\{|\ell|\le\varepsilon\}\cap(E_H\cup E_k)\mid B)\le \varepsilon n^{-30}+e^{-c_{\rm rare} n}
\]
with an absolute \(c_{\rm rare}>0\).
\end{lemma}
\begin{proof}
Fix \(B\). We use the Fourier majorant \eqref{eq:major1} for the slab with a positive exponential weight for each exceptional event. At low frequencies we estimate this weight by concentration. At high frequencies the arithmetic estimate for the normal vector gives the needed exponential decay.

\step{High frequencies under a Gaussian shift} If \(\sqrt n\le|u|\le1/\varepsilon\), Lemma~\ref{le:lattice} gives
\begin{equation}\label{eq:far}
 \dist(uv,\pi\Z^n)^2\ge2c_1n.
\end{equation}
Put \(\epsilon_{\rm sh}=c_1/10\). For \(|h|^2\le\epsilon_{\rm sh}n\), subtracting \(|h|\) from the distance in \eqref{eq:far} still leaves at least \(\sqrt{c_1 n}\). Thus
\begin{equation}\label{eq:farshift}
 \Big|\prod_\iota\cos(uv_\iota+h_\iota)\Big|\le e^{-c_1n/2}.
\end{equation}
A real exponential tilt gives hyperbolic cosines instead. At most \(\epsilon_{\rm sh}n\) coordinates have \(|h_\iota|>1\). Their total contribution to the squared distance in \eqref{eq:far} is at most \(\pi^2\epsilon_{\rm sh}n/4\). The remaining coordinates contribute at least \(c_1n\). Using \(\sin^2 x\ge4\dist(x,\pi\Z)^2/\pi^2\) and \(|\cosh(y+ix)|^2=\cosh^2 y-\sin^2 x\) gives
\begin{equation}\label{eq:cosh}
 \Big|\prod_\iota\cosh(h_\iota+iuv_\iota)\Big|
 \le e^{|h|^2/2}\exp\!\left(-\frac{\sum_{|h_\iota|\le1}\sin^2(uv_\iota)}{2\cosh^2 1}\right)
 \le e^{|h|^2/2-c_{\rm tilt}n},
\end{equation}
where \(c_{\rm tilt}=2c_1/(\pi^2\cosh^2 1)\).

\step{Fluctuation of the tail} Set \(Q_H=X^\top H X-\tr H\) and \(\lambda=n^{1.26}\).
By \eqref{eq:tails}, \(\lambda\|H\|=O(n^{-0.34})\), \(\lambda\tr H=O(n^{0.46})\), and \(\lambda^2\tr H^2=O(n^{0.12})\).
Introduce \(h\sim N(0,2\lambda H)\) independent of \(X\). Its Gaussian moment generating function and characteristic function yield, respectively,
\begin{align*}
 \E e^{\lambda X^\top H X+iu\ell}&=\E_h\prod_\iota\cosh(h_\iota+iuv_\iota),\\
 \E e^{-\lambda X^\top H X+iu\ell}&=\E_h\prod_\iota\cos(h_\iota+uv_\iota).
\end{align*}
For \(s=1/2,1\), the determinant formula and \(\lambda\|H\|=o(1)\) give
\(\E e^{s|h|^2}=\det(I-4s\lambda H)^{-1/2}\le e^{4s\lambda\tr H}\).
On \(|h|^2\le\epsilon_{\rm sh}n\), apply \eqref{eq:farshift} or \eqref{eq:cosh}. On the complement use
\[
 \E[e^{|h|^2/2}\one{|h|^2>\epsilon_{\rm sh}n}]\le e^{-\epsilon_{\rm sh}n/2}\E e^{|h|^2}.
\]
The centering factor is at most \(e^{\lambda\tr H}=e^{o(n)}\). Thus the high-frequency transforms \(\E e^{\pm\lambda Q_H+iu\ell}\) have modulus at most \(e^{-b_H n+o(n)}\), with \(b_H=\min(c_1/2,c_{\rm tilt},\epsilon_{\rm sh}/2)\).
At low frequencies their modulus is at most \(e^{Cn^{0.12}}\) by Lemma~\ref{le:HW}.
On either tail \(\pm Q_H>n^{-1.1}\), the corresponding weight \(e^{\pm\lambda Q_H}\) is at least \(e^{n^{0.16}}\). Applying \eqref{eq:major1}, each tail in the slab has probability at most
\[
 C e^{-n^{0.16}}\bigl(\varepsilon\sqrt n\,e^{C n^{0.12}}+e^{-b_H n+o(n)}\bigr).
\]
Here the low-frequency interval has length at most \(2\sqrt n\), and \(\varepsilon\) times the length of the high-frequency interval is at most two.

\step{A lower bound for the first projections} Write \(\Pi_k=\sum_{j\le k}u_ju_j^\top\). Take
\(Y=e^{-\lambda_0 X^\top \Pi_kX}\), where \(\lambda_0=\min(c_{\rm hw},1/(4C_{\rm hw}))\).
At low frequencies, Lemma~\ref{le:HW} gives \(|\E Ye^{iu\ell}|\le\E Y\le e^{-3\lambda_0 k/4}\).
At high frequencies, use \eqref{eq:farshift} with \(h\sim N(0,2\lambda_0\Pi_k)\). Its squared length is \(2\lambda_0\chi_k^2\), so
\[
 \PP(|h|^2>\epsilon_{\rm sh}n)
 \le 2^{k/2}e^{-\epsilon_{\rm sh}n/(8\lambda_0)}.
\]
Take \(b_k=\min(c_1/2,\epsilon_{\rm sh}/(8\lambda_0))\). Since \(k=o(n)\), the transform on this interval is at most \(e^{-b_k n+o(n)}\) in modulus. On \(E_k\), the weight satisfies \(Y\ge e^{-\lambda_0k/2}\). The slab contribution is therefore at most
\[
 C\varepsilon\sqrt n\,e^{-\lambda_0 k/4}+e^{-b_k n+o(n)} .
\]
The low-frequency terms for both exceptional events sum to at most \(\varepsilon n^{-30}\) because \(k=\lfloor n^\kappa\rfloor\). Their high-frequency terms sum to at most \(e^{-c_{\rm rare}n}\) with \(c_{\rm rare}=\tfrac12\min(b_H,b_k)\) for large \(n\).
\end{proof}

\subsection{Proof of Proposition~\ref{pr:flat}}
Choose an even Schwartz probability density \(\varrho\) whose Fourier transform is supported on \([-1,1]\). For example, take the normalized square of the inverse transform of a nonzero real even smooth function supported on \([-1/2,1/2]\). Let \(N_0,\ldots,N_m\) be iid with this density and independent of \(A\). Set
\[
 d=c_2/8,\qquad \vartheta=e^{-dn}.
\]
Here \(c_2=\min(c_0,c_1)\) from Corollary~\ref{co:slab}. We use the width tolerance \(\Delta\) fixed in \eqref{eq:params} and first work throughout \(e^{-dn/10}\le t\le n^{-20}\). The final lower endpoint will be chosen after we have bounded every additive error.

\step{1. Removing blocks with poor control} Retain \(\Hcal=\{B\in\Gcal:a_1\ge n^{-4}\}\). We claim
\begin{equation}\label{eq:discard}
 \E[\ind_{\{B\notin\Hcal\}}\PP(X_n\le t\mid B)]\le Ct n^{-3}+Ce^{-c_{\rm block} n}.
\end{equation}
Let \(c_{\rm rect}>0\) be a rate for the exponential error in Lemma~\ref{le:RV}(a).
The blocks outside \(\Acal\) cost at most \(e^{-c_0n}\) by Proposition~\ref{pr:LCD}. The blocks with \(a_1\le2t\) cost at most \(Ct^2+e^{-c_{\rm rect}n}\).
On the remaining blocks, all denominators in \(W_t^2\) are bounded away from zero. Since \(\sum_jy_j^2\le n\), the crossing formula \eqref{eq:cross} and Corollary~\ref{co:slab}(a) give
\begin{equation}\label{eq:condcrude}
 \PP(X_n\le t\mid B)\le Ct(1+\sqrt n/a_1)+Ce^{-c_2n}.
\end{equation}
For blocks with \(2t<a_1<n^{-4}\), integration of the rectangular lower-tail estimate gives
\[
 \E[a_1^{-1}\one{2t<a_1<n^{-4}}]
 \le n^4\PP(a_1<n^{-4})+\int_{2t}^{n^{-4}}\PP(a_1\le x)x^{-2}\,dx
 \le Cn^{-4}+Ce^{-c_{\rm rect}n}/t.
\]
Multiplication by the \(t\sqrt n\) in \eqref{eq:condcrude} yields \(O(tn^{-3.5}+\sqrt n\,e^{-c_{\rm rect}n})\).
Finally, when \(a_1\ge n^{-4}\) and \(B\in\Acal\setminus\Gcal\), the conditional bound \eqref{eq:condcrude} is at most \(Ct n^{4.5}+Ce^{-c_2n}\). Multiply this by \(\PP(B\notin\Gcal)\le n^{-100}\). Collecting the terms proves \eqref{eq:discard} with \(c_{\rm block}=\tfrac12\min(c_0,c_2,c_{\rm rect})\).

\step{2. Comparing the original and smoothed slabs} Fix \(B\in\Hcal\). Recall \(z_j=y_j/a_j\), and write \(z=(z_1,\ldots,z_m)\). Define
\[
 L=\ell+\vartheta N_0,\qquad \eta=r(N_1,\ldots,N_m),\qquad Z=z+\eta,\qquad
 W_0^*(x)=(n^{-1}+\tr H+|x|^2)^{1/2}.
\]
In \(W_0^*(z)^2\) the first \(m\) projections are kept, while \(X^\top H X\) is replaced by its mean. We compare the actual width \(W_t\) to \(W_0^*(Z)\) one change at a time. The bound \(a_1\ge n^{-4}\) gives
\begin{equation}\label{eq:smallWt}
 |z|\le n^{4.5},\quad W_0\le W_t\le M_0/4,\quad
 0\le W_t-W_0\le W_0t^2/a_1^2\le n^{-26}.
\end{equation}
On \(E_H^c\cap E_k^c\), we have \(W_0^*(z)\ge\max(w_0,|z|)\) by \eqref{eq:w0}. Dividing the difference of squared widths by the sum of the widths gives
\[
 |W_0-W_0^*(z)|\le n^{-1.1}/w_0,\qquad
 |W_0^*(Z)-W_0^*(z)|\le 2r n^{0.01}+r^2m n^{0.02}/w_0
\]
provided \(|\langle z,\eta\rangle|\le r n^{0.01}|z|\) and \(|\eta|^2\le r^2m n^{0.02}\). For each fixed \(X\), these last two noise bounds fail with conditional probability at most \(n^{-200}\) by \eqref{eq:mom}. The normal noise satisfies \(\vartheta|N_0|\le tn^{-2}\) except with probability at most \(C n^2\vartheta/t\le Ce^{-dn/2}\).
Since \(1/w_0\le n^{2\kappa}\) and \(2\kappa\le0.001\), the three width errors are at most \(n^{-26}\), \(n^{-1.099}\), and \(2n^{-1.065}+n^{-1.329}\), respectively. For large \(n\) they sum to less than \(\Delta/2\). On this event,
\begin{equation}\label{eq:bracketbounds}
 |W_t-W_0^*(Z)|\le\Delta/2,\qquad |L-\ell|\le tn^{-2},\qquad |Z|\le M_0.
\end{equation}
For the two enlarged and reduced slabs define
\[
 I_\pm=\PP\bigl(|L|\le t(W_0^*(Z)\pm\Delta),\ |Z|\le M_0\mid B\bigr).
\]
Both widths are positive since \(W_0^*(x)\ge n^{-1/2}>\Delta\).
The bounds \eqref{eq:bracketbounds} and \(n^{-2}\le\Delta/2\) give the two inclusions. Indeed, a crossing gives
\(|L|\le tW_t+tn^{-2}\le t(W_0^*(Z)+\Delta)\). Conversely the reduced slab gives
\(|\ell|\le t(W_0^*(Z)-\Delta)+tn^{-2}\le tW_t\).

To bound errors in these inclusions, we need only bound failures on a crossing or on the reduced slab. A crossing forces \(|\ell|\le2M_0t\) by \eqref{eq:smallWt}. The reduced slab does as well whenever the normal noise satisfies its bound, since \(|Z|\le M_0\) there. Thus Lemma~\ref{le:rare} applies at \(\varepsilon=2M_0t\in[e^{-c_0n},1]\). To handle the tangential noises, condition first on \(X\) to fix the slab indicator. For every such \(X\), the conditional failure probability of these bounds is at most \(n^{-200}\). Averaging over \(X\) gives a bound \(Cn^{-200}(M_0t+e^{-c_2n})\) by Corollary~\ref{co:slab}(a). Including the normal-noise failure gives
\begin{equation}\label{eq:Ipm}
 I_- - e_B\le\PP(X_n\le t\mid B)\le I_++e_B,\qquad
 e_B\le Ct n^{-20}+C\big(e^{-c_{\rm rare}n}+e^{-dn/2}+e^{-c_2n}\big).
\end{equation}

\step{3. Evaluating the density at the center of the slab} Let \(p_B^{\rm flat}(s,x)\) be the joint density of \((L,Z)\) for this smoothing, conditional on \(B\). Then
\begin{equation}\label{eq:slabint}
 I_\pm=\int_{|x|\le M_0}\int_{-b_\pm(x)}^{b_\pm(x)}p_B^{\rm flat}(s,x)\,ds\,dx,
 \qquad b_\pm(x)=t(W_0^*(x)\pm\Delta)\le2M_0t .
\end{equation}
We will replace the integrand by its value at \(s=0\). Direct differentiation would see the very narrow smoothing width \(\vartheta\). Truncating frequencies in the normal coordinate gives a bound uniform in that width.

Use the Fourier convention \(\widehat p(u,\nu)=\int e^{i(us+\langle\nu,x\rangle)}p(s,x)\,ds\,dx\). Independence of the noises gives
\[
 \widehat{p_B^{\rm flat}}(u,\nu)=\varphi(w)\widehat\varrho(\vartheta u)\prod_{j\le m}\widehat\varrho(r\nu_j),
 \qquad w=uv+\sum_{j\le m}(\nu_j/a_j)u_j .
\]
On the Fourier support, \(|u|\le e^{dn}\) and \(|\nu_j|\le1/r\). Consequently \(|w|\le e^{dn}+n^6\) and \(|Bw|=|\nu|/\sqrt n\le\sqrt m/(r\sqrt n)\le n^{0.975}\).
For \(|u|\ge\sqrt n\) we also have \(|w|\ge\sqrt n\), so all three inequalities in \eqref{eq:latticecond} hold for large \(n\). This yields \(|\varphi(w)|\le e^{-c_1n}\).

Take an even smooth \(\chi\) with values in \([0,1]\), equal to one on \([-1,1]\) and zero off \([-2,2]\). Multiply \(\widehat{p_B^{\rm flat}}\) by \(\chi(u/\sqrt n)\), and denote the inverse transform by \(p_<\). Fourier inversion and \(m\log n=o(n)\), \(d\le c_1/8\), give
\[
 \|p_B^{\rm flat}-p_<\|_\infty
 \le (2\pi)^{-m-1}e^{-c_1n}(2e^{dn})(2/r)^m\le e^{-c_1n/2}.
\]
The integration volume in \eqref{eq:slabint} is at most \(t(CM_0)^{m+1}\), so discarding and then restoring this high-frequency part costs \(Ct e^{-c_1n/4}\).
For the remaining frequencies, the cutoff convolves the density in its first coordinate with \(\sqrt n\,\check\chi(\sqrt n\,\cdot)\), where \(\check\chi=\mathcal F^{-1}\chi\). Its derivative is bounded by \(n\|\check\chi'\|_\infty\). Integrating the original joint density in the first coordinate gives the marginal \(p_Z(x)\). Consequently
\[
 |p_<(s,x)-p_<(0,x)|\le C n |s|p_Z(x).
\]
Integrating this over \eqref{eq:slabint} costs at most \(CnM_0^2t^2\).
Finally we can bound the marginal density of \(L\) at zero by Fourier inversion. Corollary~\ref{co:slab} bounds the integrand by \(e^{-c_1u^2}\) up to \(\sqrt n\), and by \(e^{-c_1n}\) from there to \(e^{dn}\). Outside this range \(\widehat\varrho(\vartheta u)=0\). This proves \(p_L(0)\le C\). Changing either half-width by \(t\Delta\) costs at most \(2t\Delta p_L(0)\). In particular,
\begin{equation}\label{eq:flatcoef}
 I_\pm=2tC_B+O(t\Delta+nM_0^2t^2+t e^{-c_1n/4}),\qquad
 C_B=\int_{|x|\le M_0}W_0^*(x)p_B^{\rm flat}(0,x)\,dx.
\end{equation}
The coefficient \(C_B\) is nonnegative, finite, and independent of \(t\).

\step{4. Averaging and choosing the lower endpoint} Set
\(a_n=2\E[\ind_{\{B\in\Hcal\}}C_B]\), assigning the integrand zero on the discarded blocks. The following bounds are in absolute probability units, up to absolute factors. After division by \(t\), only the last column can grow as the threshold decreases.
\begin{center}
\small
\begin{tabular}{@{}p{.37\linewidth}p{.25\linewidth}p{.29\linewidth}@{}}
\hline
\textbf{Source} & \textbf{Error with a factor \(t\)} & \textbf{Additive exponential error}\\
\hline
Discarded blocks, \eqref{eq:discard} & \(tn^{-3}\) & \(e^{-c_{\rm block}n}\)\\
\(E_H\cup E_k\) on the slab & \(2M_0t n^{-30}\) & \(e^{-c_{\rm rare}n}\)\\
Tangential noise on the slab & \(M_0t n^{-200}\) & \(n^{-200}e^{-c_2n}\)\\
Normal noise & \(0\) & \(e^{-dn/2}\)\\
Width tolerance in \eqref{eq:flatcoef} & \(t\Delta\) & \(0\)\\
Variation of the density & \(nM_0^2t^2\) & \(0\)\\
High-frequency cutoff & \(t e^{-c_1n/4}\) & \(0\)\\
\hline
\end{tabular}
\end{center}
In the middle column, \(M_0=n^6\), \(nM_0^2t\le n^{-7}\), and \(\Delta=n^{-1.05}\), so the sum is \(O(t n^{-1.05})\). In the last column take
\[
 c_{\rm add}=\min(c_{\rm block},c_{\rm rare},c_2,d/2),\qquad
 \cflat=\min(d/10,c_{\rm add}/2).
\]
Equations \eqref{eq:discard}, \eqref{eq:Ipm}, and \eqref{eq:flatcoef} now give
\[
 |F_n(t)-a_nt|\le Ct n^{-1.05}+Ce^{-c_{\rm add}n}\qquad(e^{-dn/10}\le t\le n^{-20}).
\]
For \(t\ge e^{-\cflat n}\) the last term divided by \(t\) is at most \(Ce^{-c_{\rm add}n/2}\). This proves the claimed uniform estimate. At \(t=n^{-20}\) the square-matrix bound \(F_n(t)\le Ct+e^{-cn}\) from \cite{rvsmall} also gives \(a_n\le C\).
\section{Polynomial scales}\label{sec:poly}
Section~\ref{sec:flat} approximates \(F_n(t_*)/t_*\) by \(a_n\) at \(t_*=n^{-20}\). We now compare this ratio with \(F_n(t)/t\) for \(t_*\le t\le n^{-\theta}\). By row exchangeability we may count a crossing with weight \(U^2=n|u_n|^2\), where \(u\) is the least left singular vector sampled as in Lemma~\ref{le:exchange}. This weight serves two purposes. It limits the contribution from blocks with a small \(a_1\), and it cancels a factor in the derivative of the crossing map \(s\mapsto sW_s\). Recall the absolute tolerance \(\Delta\) and relative tolerance \(\epsilon_n\) from \eqref{eq:params}. Set
\[
 T=t\log^3 n,\qquad
 \Hcal_t=\{B\in\Gcal:a_1\ge\max(T,n^{-4})\}.
\]
On \(\Hcal_t\), every denominator \(1-s^2/a_j^2\) is bounded away from zero for \(0\le s\le 2t\). At the anchor \(t_*\), the cutoff \(T\) is smaller than \(n^{-4}\), so \(\Hcal_{t_*}=\Hcal\). For the next proposition the same number \(\bar D_n\) must work at every scale. The lower bound at \(t_*\) will let us use the value obtained in Corollary~\ref{co:slope}.
\begin{proposition}\label{pr:poly}
For each large \(n\) there is \(\bar D_n\ge0\) such that
\begin{enumerate}[label=\textup{(\alph*)}]
\item \(F_n(t)\le(1+3\epsilon_n)\bar D_n t+t n^{-6}\) for \(t_*\le t\le n^{-\theta}\).
\item \(F_n(t_*)\ge(1-3\epsilon_n)\bar D_n t_*-t_* n^{-6}\).
\end{enumerate}
\end{proposition}
\begin{corollary}\label{co:poly}
For \(t_*\le t\le n^{-\theta}\), one has \(F_n(t)\le t(a_n+Cn^{-1.02})\le t(1-\frac1{2n}+Cn^{-1.02})\).
\end{corollary}
\begin{proof}
At \(t_*\) the weighted lower bound of Proposition~\ref{pr:poly}(b) and the unweighted upper bound of Proposition~\ref{pr:flat} concern the same probability. Dividing both by \(t_*\) gives
\[
 (1-3\epsilon_n)\bar D_n
 \le F_n(t_*)/t_*+n^{-6}
 \le a_n+Cn^{-1.05}+n^{-6}.
\]
Since \(a_n\le C\), division by \(1-3\epsilon_n\) gives \(\bar D_n\le a_n+Cn^{-1.02}\). Proposition~\ref{pr:poly}(a) uses this same \(\bar D_n\) for all thresholds in the polynomial range. Substitution and Corollary~\ref{co:slope} prove the claim.
\end{proof}

\subsection{Blocks with a small singular value}
\begin{lemma}\label{le:smallblock}
For \(t_*\le t\le n^{-\theta}\),
\[
 \E[U^2\one{X_n\le t}\one{B\notin\Hcal_t}]\le t n^{-7}+t^{10/3}n^{3\kappa}
\]
for large \(n\).
\end{lemma}
\begin{proof}
Set \(Y_t=U^2\one{X_n\le t}\). Since \(Y_t\le n\), blocks outside \(\Gcal\) contribute at most \(n^{-99}\). Fix a block in \(\Gcal\) for the conditional estimates below. By \eqref{eq:crudeweight},
\[
 Y_t\le\min\{n,W_0^{-2},t^2/\ell^2\},
\]
where the last term is interpreted as infinity at \(\ell=0\). For any event \(E\) and deterministic cutoff \(M>0\), integration of level sets gives
\[
 \E[\one{E}\min(M,t^2/\ell^2)\mid B]
 =\int_0^M\PP(E,\ |\ell|\le t/\sqrt h\mid B)\,dh.
\]
For \(E\) equal to the whole space, Corollary~\ref{co:slab}(a) bounds this by \(Ct\sqrt M+CM e^{-c_2 n}\). We choose \(M\) according to the size of the block and its first tangential projection. Cases (ii) and (iii) concern \(n^{-4}\le a_1<T\). They are empty when \(T\le n^{-4}\). In both cases set \(\beta=T^{1/3}\).

\begin{enumerate}[label=\textup{(\roman*)}]
\item \emph{A very small block singular value, \(a_1<n^{-4}\).}
With \(M=n\), the conditional bound is \(Ct\sqrt n+Cn e^{-c_2 n}\).
Lemma~\ref{le:RV}(a) gives \(\PP(a_1<n^{-4})\le Cn^{-8}\), so the averaged contribution is \(O(tn^{-7.5})\).

\item \emph{An intermediate block with \(|y_1|>\beta\).}
Here \(a_1<T\le\beta\le1\) for large \(n\). On \(|y_1|>\beta\) we have \(W_0\ge |y_1|/a_1\ge\beta/a_1\). With \(M=(a_1/\beta)^2\le1\), the conditional bound becomes
\[
 \E[Y_t\one{|y_1|>\beta}\mid B]
 \le Ct a_1/\beta + C e^{-c_2 n}
 \le CtT^{2/3}+C e^{-c_2 n}.
\]

\item \emph{An intermediate block with \(|y_1|\le\beta\).}
On \(E_k^c\), \eqref{eq:w0} gives the cutoff \(M=w_0^{-2}\). The part of the level integral with \(h<t^2\) is at most \(t^2\). On \([t^2,M]\), the width \(t/\sqrt h\) lies in \([tw_0,1]\subset[n^{-100},1]\). Also \(\beta\ge a_1\ge n^{-4}\). Thus Corollary~\ref{co:slab}(b) bounds the integrand by \(Ct\beta/\sqrt h+Ce^{-c_2 n}\). On \(E_k\), use \(Y_t\le n\) and \(\PP(E_k\mid B)\le2e^{-ck}\). Together these give
\[
 \E[Y_t\one{|y_1|\le\beta}\mid B]
 \le t^2+Ct\beta/w_0+C e^{-c_2n/2}+2ne^{-ck},
\]
where \(w_0^{-2}\le n^{4\kappa}\) absorbs only a polynomial factor into the exponential error.
\end{enumerate}

To average (ii) and (iii), use \(\PP(a_1<T)\le CT^2+e^{-cn}\) and \(w_0^{-1}\le n^{2\kappa}\). Their algebraic contributions are bounded by
\[
 C\big(tT^{8/3}+t^2T^2+tT^{7/3}n^{2\kappa}\big).
\]
Since \(t\le T\le1\), the last term dominates. Using \(T=t\log^3 n\), its ratio to \(t^{10/3}n^{3\kappa}\) is \(n^{-\kappa}\log^7n\to0\). All the exponential errors are \(o(tn^{-7})\) because \(t\ge n^{-20}\). Adding (i) and the blocks outside \(\Gcal\) proves the claim.
\end{proof}

\subsection{A smoothed crossing process}
For each retained block, \eqref{eq:cross} identifies the first crossing of \(sW_s=|\ell|\) with the rescaled least singular value, as long as the crossing occurs below \(a_1\). To integrate over the possible crossing points, we add independent noise to the normal coordinate and the first \(m\) tangential coordinates. The contribution of the remaining tangential coordinates at zero is \(X^\top H X\). Its conditional mean is \(\tr H\).
Take fresh noise \(N_1,\ldots,N_m\) with density \(\varrho\) and a standard Gaussian \(g_0\), independent of one another and \(A\). Write \(\sigma_N^2=\E N_1^2\) for the noise variance. Set \(\eta=r(N_j)_{j\le m}\) and \(z=(y_j/a_j)_{j\le m}\). On \(B\in\Hcal\), put
\[
 L=\ell+n^{-90}g_0,\quad Z=z+\eta,\quad b=n^{-1}+\tr H,\quad
 \lambda_j(s)=(1-s^2/a_j^2)^{-1}.
\]
For \(0\le s\le a_1/2\) and \(x\in\R^m\), define
\[
 W_s^*(x)^2=b+\sum_{j\le m}\lambda_j(s)x_j^2,\qquad
 V_s^*(x)^2=b+\sum_{j\le m}\lambda_j(s)^2x_j^2 .
\]
The derivative of \(sW_s^*(x)\) is \(V_s^*(x)^2/W_s^*(x)\ge W_s^*(x)>0\). Write \(p_B^{\rm pol}\) for the conditional density of \((L,Z)\) on this block. The superscript distinguishes the polynomial-range smoothing from the one in Section~\ref{sec:flat}. When \(B\in\Hcal_t\), \(\lambda_j(s)\le1+5\log^{-6}n\) for \(0\le s\le2t\). For \(|l|\le2t W_{2t}^*(x)\), let \(s^*(l,x)\) be the unique solution in \([0,2t]\) of \(sW_s^*(x)=|l|\), and set \(s^*=\infty\) otherwise. For \(0\le h\le2t\), write \(G_B(h)=\E[V_{s^*}^*(Z)^{-2}\one{s^*\le h}\mid B]\), with integrand defined to be zero when \(s^*=\infty\). Define
\begin{equation}\label{eq:crossdensity}
 D_B(s)=\int_{\R^m}
 \frac{p_B^{\rm pol}(sW_s^*(x),x)+p_B^{\rm pol}(-sW_s^*(x),x)}{W_s^*(x)}\,dx
 \qquad(0\le s\le a_1/2).
\end{equation}
\begin{lemma}\label{le:crossdensity}
For \(B\in\Hcal_t\), \(G_B(h)=\int_0^hD_B(s)\,ds\) on \(0\le h\le2t\). Moreover,
\begin{equation}\label{eq:DFourier}
 D_B(s)=\frac1\pi\int_{\R}\operatorname{Re}\E\left[
 e^{iuL}\frac{\cos(usW_s^*(Z))}{W_s^*(Z)}\,\middle|\,B\right]du
 \qquad(0\le s\le a_1/2).
\end{equation}
The integral is absolutely convergent.
\end{lemma}
\begin{proof}
For fixed \(Z=x\), split the integral in the normal coordinate into \(l\ge0\) and \(l<0\). On each half, substitute \(|l|=sW_s^*(x)\). The cancellation that motivates the choice of weight is
\[
 \frac{1}{V_s^*(x)^2}\frac{d}{ds}\big[sW_s^*(x)\big]
   = \frac{1}{V_s^*(x)^2}\frac{V_s^*(x)^2}{W_s^*(x)}
   = \frac{1}{W_s^*(x)}.
\]
Integrating the two contributions against the joint density gives \(G_B(h)=\int_0^hD_B(s)\,ds\) by Tonelli's theorem. Let \(p_Z\) be the marginal density of \(Z\) on this block. For the Fourier representation, the partial Fourier transform of \(p_B^{\rm pol}(l,x)\) in \(l\) is bounded in modulus by \(e^{-n^{-180}u^2/2}p_Z(x)\) because \(g_0\) is independent of \((\ell,Z)\). Invert this transform at \(l=\pm sW_s^*(x)\) and use \(W_s^*(x)\ge\sqrt b\) to apply Fubini's theorem.
\end{proof}
\begin{lemma}[Comparison]\label{le:compare}
For \(B\in\Hcal_t\) and \(t_*\le t\le n^{-\theta}\),
\[
 (1-\epsilon_n)G_B((1-\epsilon_n)t)-n^{-90}
 \le \E[U^2\one{X_n\le t}\mid B]
 \le (1+\epsilon_n)G_B((1+\epsilon_n)t)+n^{-90}.
\]
\end{lemma}
\begin{proof}
\step{Control of the perturbation} Condition on \(B\in\Hcal_t\). By Lemma~\ref{le:HW} and \eqref{eq:tails}, \(\PP(E_H\cup E_k)\le n^{-200}\) for large \(n\). Impose additionally
\[
 n^{-90}|g_0|\le n^{-80},\qquad |\eta|^2\le r^2m n^{0.02},\qquad
 \sup_{0\le s\le2t}\Big|\sum_{j\le m}\lambda_j(s)^i z_j\eta_j\Big|
       \le2r n^{0.01}|z| \quad (i=1,2).
\]
Given \(X\), put \(M_j^{\eta}=\sum_{h\le j}z_h\eta_h\). For each \(s\), the coefficients \(\lambda_j(s)^i\) decrease in \(j\). Summation by parts gives, simultaneously for all \(0\le s\le2t\) and \(i=1,2\),
\[
 \left|\sum_{j\le m}\lambda_j(s)^i z_j\eta_j\right|
 \le (1-s^2/a_1^2)^{-i}\max_{j\le m}|M_j^{\eta}|
 \le 2\max_{j\le m}|M_j^{\eta}|.
\]
By \eqref{eq:mom}, \((\E[|M_j^{\eta}|^p\mid X,B])^{1/p}\le C_p r|z|\). Markov's inequality and a union bound over \(m\le n\) give
\(\PP(\max_j|M_j^{\eta}|>r n^{0.01}|z|\mid X,B)\le C_p n^{1-0.01p}\) when \(z\ne0\). All sums vanish for \(z=0\). Taking \(p=21000\), and using \eqref{eq:mom} for \(|\eta|^2\) and the Gaussian tail for \(g_0\), bounds the total conditional failure probability by \(n^{-200}\) for large \(n\).

On the joint good event, \(W_s^2-W_s^*(Z)^2\) equals
\[
 X^\top H X-\tr H+\sum_{j>m}(\lambda_j(s)-1)z_j^2
             -\sum_{j\le m}\lambda_j(s)(2z_j\eta_j+\eta_j^2).
\]
The first term is bounded by \(n^{-1.1}\), the second by \(Ct^2m^{-2}X^\top H X\le n^{-2}\), and the third in absolute value by \(4r n^{0.01}|z|+2r^2m n^{0.02}\). The same argument applies to the \(V\)'s with \(\lambda_j(s)^2\). To pass from squares to square roots, divide the term proportional to \(|z|\) by the lower bound \(W_s,V_s\ge |z|\), and use \(W_s,V_s\ge w_0\) for all the other terms. Each absolute error is at most
\[
 C\left[(n^{-1.1}+n^{-2}+r^2m n^{0.02})/w_0 + r n^{0.01}\right]\le \Delta,
\]
since \(w_0^{-1}\le n^{2\kappa}\), \(r n^{0.01}=n^{-1.065}\), and \(r^2m n^{0.02}\le n^{-1.33}\). Consequently
\begin{equation}\label{eq:WVclose}
 \sup_{0\le s\le2t}\bigl(|W_s-W_s^*(Z)|+|V_s-V_s^*(Z)|\bigr)\le2\Delta,
 \qquad |L-\ell|\le n^{-80}.
\end{equation}
\step{Comparison of the crossing points}
Write \(g(s)=sW_s\) and \(g^*(s)=sW_s^*(Z)\). On the event just constructed, \(W_s,V_s\ge w_0\) and \eqref{eq:WVclose} gives \(W_s^*(Z),V_s^*(Z)\ge w_0/2\). Thus on \([0,2t]\),
\[
 |g-g^*|\le4t\Delta,\qquad g'\ge w_0/2,\qquad (g^*)'\ge w_0/2 .
\]
A vertical change of size \(e_t=4t\Delta+n^{-80}\) moves a crossing horizontally by at most \(d_t=2e_t/w_0\). Since \(w_0\ge n^{-2\kappa}\), \(\Delta=n^{-1.05}\), and \(t\ge n^{-20}\), we have \(d_t\le t\epsilon_n/4\) for large \(n\).

Let \(s_0=X_n\). If \(s_0\le t\), then \(g(s_0)=|\ell|\) by \eqref{eq:cross} and continuity. The perturbation bound gives \(\big||L|-g^*(s_0)\big|\le e_t\), so
\[
 g^*((s_0-d_t)_+)\le |L|\le g^*(s_0+d_t),\qquad s_0+d_t\le(1+\epsilon_n)t<2t.
\]
When \(s_0<d_t\), the lower inequality uses \(g^*(0)=0\). Otherwise both inequalities follow by integrating the derivative bound. There is therefore a smoothed crossing \(s^*=s^*(L,Z)\) and \(|s^*-s_0|\le d_t\).
Conversely, suppose \(s^*\le(1-\epsilon_n)t\). From \(g^*(s^*)=|L|\) we obtain \(\big||\ell|-g(s^*)\big|\le e_t\) and
\[
 g((s^*-d_t)_+)\le |\ell|\le g(s^*+d_t),\qquad s^*+d_t\le t<a_1.
\]
The unique root of \(g(s)=|\ell|\) in this interval equals \(X_n\) by \eqref{eq:cross}. Again \(|s_0-s^*|\le d_t\).

\step{Comparison of the weights}
At the same argument \(s=s_0\), \eqref{eq:WVclose} and the lower bound on both square roots imply
\[
 |\log V_{s_0}-\log V_{s_0}^*(Z)|\le4\Delta/w_0=o(\epsilon_n).
\]
For the change of argument, direct differentiation of \(V_s^*(x)^2\) gives
\[
 0\le\partial_s\log V_s^*(x)
 \le\frac{2s}{a_1^2(1-s^2/a_1^2)}\le\frac{5t}{a_1^2}\quad(0\le s\le2t).
\]
Indeed each logarithmic derivative of \(\lambda_j(s)\) is at most the middle quantity. Since \(a_1\ge t\log^3 n\), the change from \(s_0\) to \(s^*\) costs at most \(5td_t/a_1^2=o(\epsilon_n)\). By \eqref{eq:weight}, \(U^2=V_{s_0}^{-2}\) whenever \(s_0\le t<a_1\). Exponentiating the sum of the two log bounds yields
\[
 (1-\epsilon_n)V_{s^*}^*(Z)^{-2}\le U^2\le(1+\epsilon_n)V_{s^*}^*(Z)^{-2}
\]
whenever either comparison of crossing points above applies. Combine this with the two event inclusions. On the exceptional event each weight is at most \(n\), so taking conditional expectations costs at most \(C n^{-199}\le n^{-90}\) on each side.
\end{proof}

\subsection{High frequencies}
To obtain the curvature of \(D_B(s)\), we will expand the cosine in \eqref{eq:DFourier} for \(|u|\le K\). The next lemma justifies that truncation. Its difficulty is that the weight \(f_u(q_s)\) and the oscillating factor \(e^{iu\ell}\) both depend on the same sign row \(X\). To apply the arithmetic estimate of Lemma~\ref{le:lattice}, we will express their product as a sum of integrals of \(\E[e^{i\langle w,X\rangle}X^S]\), where \(X^S\) is a monomial of bounded degree.
\begin{lemma}\label{le:high}
For \(B\in\Hcal_t,\ t_*\le t\le n^{-\theta},\ 0\le s\le2t\), put \(q_s=W_s^*(Z)^2\). For each of the functions \(f_u(q)=q^{-1/2}\) and \(f_u(q)=q^{-1/2}\cos(us\sqrt q)\),
\[
 \int_{|u|>K}\big|\E[e^{iuL}f_u(q_s)\mid B]\big|\,du\le n^{-100}
\]
for sufficiently large \(n\).
\end{lemma}
\begin{proof}
Fix \(B\in\Hcal_t\) and \(s\in[0,2t]\). Write \(\lambda_j=\lambda_j(s)\in[1,2]\). All expectations in the proof are conditional on this fixed block. Independence of \(g_0\) gives
\(\E[e^{iuL}f_u(q_s)]=e^{-n^{-180}u^2/2}\E[e^{iu\ell}f_u(q_s)]\).
Since \(|f_u(q_s)|\le\sqrt n\), the integral over \(|u|\ge n^{91}\) is at most \(e^{-n}\). Fix \(K<|u|<n^{91}\).

\step{1. Bound derivatives of the weight} The weight \(q^{-1/2}\) is large near zero. We first show that the first few tangential coordinates already keep \(q_s\) away from zero except on an event of negligible probability. Set \(k'=\lceil\log^2 n\rceil,\ q_{\rm cut}=\log^{-9}n\). Except on an event of probability \(n^{-300}\), \(\sum_{j\le k'}y_j^2\ge k'/2\) by Lemma~\ref{le:HW}, and \(\sum_{j\le k'}N_j^2\le n^{0.02}\) by \eqref{eq:mom}. Then (G3) gives
\[
 \sqrt{q_s}\ge\Big(\sum_{j\le k'}Z_j^2\Big)^{1/2}
 \ge\frac{\sqrt{k'/2}}{a_{k'}}-r n^{0.01}
 \ge c\log^{-4}n\ge\sqrt{2q_{\rm cut}}.
\]
Let \(\varpi:\R\to[q_{\rm cut},\infty)\) be smooth and nondecreasing, equal to the identity for \(q\ge2q_{\rm cut}\), with \(\|\varpi^{(i)}\|_\infty\le C_iq_{\rm cut}^{1-i}\) for \(i\ge1\). For \(\widetilde f=f_u\circ\varpi\), replacement of \(f_u(q_s)\) by \(\widetilde f(q_s)\) costs at most \(2\sqrt n\,n^{-300}\) in expectation. Repeated differentiation gives, for \(i\ge0\),
\begin{equation}\label{eq:fder}
 \|\widetilde f^{(i)}\|_\infty\le C_i q_{\rm cut}^{-1/2-i}(1+|u|)^i.
\end{equation}
Indeed differentiating \(\cos(us\sqrt q)\) introduces factors bounded using \(s\le1,\ q\ge q_{\rm cut}\), and the cutoff derivatives have the displayed scaling.

\step{2. Keep the first coordinates and center the remaining ones} The number of coordinates to keep grows with \(|u|\). The split must make \(|u|j_0^{-3/2}\) small for the Taylor expansion and keep \(j_0\log n\) below the arithmetic decay exponent \(u^2\) when \(u^2\le n\). The power \(3/2\) meets both requirements. Take \(j_0=\min(m,\lceil|u|^{3/2}\rceil)\), so \(j_0\ge\max(k',J_0)\) for large \(n\). Write \(q_s=\bar q+\zeta\), with
\begin{align*}
 \bar q&=b+\sum_{j\le j_0}\lambda_jZ_j^2+\sum_{j_0<j\le m}\lambda_j(a_j^{-2}+r^2\sigma_N^2),\\
 \zeta&=\sum_{j_0<j\le m}\lambda_j\left(\frac{y_j^2-1}{a_j^2}
                +\frac{2r y_jN_j}{a_j}+r^2(N_j^2-\sigma_N^2)\right).
\end{align*}
The sign row has identity covariance, so \(\E[y_j^2\mid B]=1\). The noise is independent of the sign row, has mean zero and has variance \(\sigma_N^2\). Thus \(\E[\zeta\mid B]=0\). For \(j_0=m\) the tail is empty, so \(\zeta=0\). For \(j_0<m\) and each fixed \(p\ge2\),
\begin{equation}\label{eq:R}
 \|\zeta\|_{L^p}\le C_p j_0^{-3/2}.
\end{equation}
For the sign quadratic form, the squared Hilbert--Schmidt norm is \(\sum_{j_0<j\le m}\lambda_j^2a_j^{-4}\le Cj_0^{-3}\) by (G2), so use Lemma~\ref{le:HW}. For the mixed term, condition on \(X\) in \eqref{eq:mom} and then use
\(\|\sum_{j_0<j\le m}y_j^2/a_j^2\|_{L^{p/2}}\le C_p/j_0\).
Since \(r\le m^{-1}\), the result is at most \(C_p j_0^{-3/2}\). The noise quadratic term is bounded by \(C_p r^2\sqrt m\le C_p m^{-3/2}\).

\step{3. Bound the Taylor remainder} Fix the expansion order \(p_{\rm T}=\lceil100/\tau\rceil\), independent of \(n\). Expand \(\widetilde f(q_s)\) about \(\bar q\) to degree \(p_{\rm T}-1\). If \(j_0<m\), \eqref{eq:fder}, \eqref{eq:R} and \(j_0\ge |u|^{3/2}\) give the factor \((1+|u|)^{p_{\rm T}}j_0^{-3p_{\rm T}/2}\le C_{p_{\rm T}}|u|^{-5p_{\rm T}/4}\). Thus the expected absolute remainder is at most
\(C_{p_{\rm T}}\log^{9p_{\rm T}+5}n\,|u|^{-5p_{\rm T}/4}\). Integration above \(K\) gives at most \(C_{p_{\rm T}}(\log n)^{9p_{\rm T}+5}K^{1-5p_{\rm T}/4}\le n^{-110}\) for large \(n\), since \(\tau p_{\rm T}\ge100\). For \(j_0=m\) the expansion is exact at degree zero.

\step{4. Integrate out the tail noise at fixed sign row}
For \(0\le i<p_{\rm T}\) put \(\mathcal T_i(u)=\E[e^{iuL}\widetilde f^{(i)}(\bar q)\zeta^i]\). The \(i\)-th Taylor term is \(\mathcal T_i(u)/i!\).
After fixing the sign row as well, the conditional mean becomes
\[
 \E[\zeta\mid X,B]=\sum_{j_0<j\le m}\lambda_j(y_j^2-1)/a_j^2.
\]
This quantity can be nonzero, so we retain the full conditional moments. Given \((X,B)\), the only remaining randomness in \(\bar q\), \(\zeta\), and \(L\) comes respectively from the noises \((N_j)_{j\le j_0}\), \((N_j)_{j_0<j\le m}\), and \(g_0\). These three families are independent. Let
\[
 \mathcal M_i(X)=\E[\zeta^i\mid X,B]
       =\sum_{|S|\le2i}\pi_S X^S,\qquad
 X^S=\prod_{\iota\in S}X_\iota,\qquad \sum_S|\pi_S|\le C_i n^{2i}.
\]
For \(i=0\), \(\mathcal M_0=1\). For every \(i\), integrating the tail noises first, and then \(g_0\), gives the identity
\[
 \mathcal T_i(u)=e^{-n^{-180}u^2/2}
     \E_{X,N_{\le j_0}}\!\left[e^{iu\ell}\widetilde f^{(i)}(\bar q) \mathcal M_i(X)\right].
\]
To bound the coefficient sum of \(\mathcal M_i\), expand each \(y_j=\sum_\iota u_{j,\iota}X_\iota\) and use \(\sum_\iota|u_{j,\iota}|\le\sqrt n\). For fixed tail noises, the sum of absolute polynomial coefficients of \(\zeta\) is bounded by
\(C\sum_{j_0<j\le m}(n+\sqrt n\,|N_j|+N_j^2+1)\), since the scalar weights are bounded on this range. Its \(L^i\) norm is at most \(C_i n^2\) for \(i\ge1\). The bounds multiply under products, and reducing \(X_\iota^2=1\) can only decrease them.

\step{5. Integrate out the first noises by Fourier inversion}
Write \(\widetilde f^{(i)}(\bar q)=g_i(Z_{\le j_0})\), with \(\|g_i\|_\infty\) bounded by \eqref{eq:fder}. Multiply \(g_i\) by a smooth cutoff equal to one on the ball of radius \(n^8\) and zero outside the ball of radius \(2n^8\), and denote the product by \(\psi_i\). Since \(|z|\le\sqrt n/a_1\le n^{4.5}\), the cutoff changes the factor only when \(|\eta|>n^8/2\), an event of probability at most \(n^{-D'}\) for any fixed \(D'\) by \eqref{eq:mom}. After the tail noises are integrated out, the polynomial in \(X\) has modulus at most \(C_i n^{2i}\). Thus the cutoff costs at most \(C_i n^{2i} q_{\rm cut}^{-1/2-i}(1+n^{91})^i n^{-D'}\). Choosing \(D'=200p_{\rm T}+300\) makes this at most \(n^{-200}\), uniformly on the frequency interval.

For each fixed \(X\), Fourier inversion of \(\psi_i\) involves the characteristic function of the first \(j_0\) noises. It contributes the product \(\prod_{j\le j_0}\widehat\varrho(r\nu_j)\), which has modulus at most one and vanishes unless \(|\nu_j|\le1/r\) for every \(j\). The sign-dependent part of the Fourier phase is linear, since \(z_j=\langle u_j,X\rangle/a_j\). Now average over \(X\). For each monomial \(X^S\), after possibly reversing \(\nu\), the resulting integral has absolute value at most
\[
 (2\pi)^{-j_0}\int_{[-1/r,1/r]^{j_0}}|\widehat\psi_i(\nu)|
       \left|\E_X[e^{i\langle w,X\rangle}X^S]\right|\,d\nu,
 \qquad w=uv+\sum_{j\le j_0}(\nu_j/a_j)u_j .
\]
The compact support gives \(|\widehat\psi_i|\le\|g_i\|_\infty(Cn^8)^{j_0}\). To apply Lemma~\ref{le:lattice}, use orthogonality and (G4) to obtain
\[
 |u|\le|w|\le n^{92},\qquad |Bw|=|\nu|/\sqrt n\le n^{0.975},
 \qquad |Bw|\le d_m|w|\le c_0\sqrt n\,|w|/2.
\]
Thus \eqref{eq:latticecond} holds for large \(n\), and \eqref{eq:walsh} gives the factor \(e^{2|S|}e^{-c_1\min(u^2,n)}\). Since \(|S|\le2i\), absorb \(e^{2|S|}\) into \(C_i\).

\step{6. Compare the Fourier volume with the arithmetic decay}
The coefficient sum of \(\mathcal M_i\) costs \(C_i n^{2i}\). The derivative bound supplies \(q_{\rm cut}^{-1/2-i}(1+|u|)^i\). The support volume of \(\psi_i\) and the volume of the frequency cube together cost \((Cn^8/r)^{j_0}\). Putting these into one bound gives
\[
 |\mathcal T_i(u)|\le C_i n^{2i}q_{\rm cut}^{-1/2-i}(1+|u|)^i(Cn^8/r)^{j_0}
       e^{-c_1\min(u^2,n)}+n^{-200}.
\]
For \(u^2\le n\), \(j_0\le2|u|^{3/2}\), so \(j_0\log(Cn^8/r)/u^2\le C\log n/\sqrt K=o(1)\). For \(u^2>n\), \(j_0\log(Cn^8/r)\le Cm\log n=o(n)\). Thus the dimensional factor costs at most half the negative exponent. The remaining polynomial in \(n\) has degree fixed once \(p_{\rm T}\) is chosen. Its product with \(e^{-c_1K^2/2}\) is \(o(n^{-200})\), so \(|\mathcal T_i(u)|\le2n^{-200}\). Integrating over an interval of length at most \(2n^{91}\) and summing the fixed number of terms costs \(O(n^{-109})\) with a constant depending only on \(p_{\rm T}\). The integrated regularisation error is \(O(n^{-208})\) and the Taylor remainder is at most \(n^{-110}\). Their sum and the Gaussian tail are at most \(n^{-100}\) for sufficiently large \(n\).
\end{proof}
\subsection{Low frequencies and the curvature term}
The bounded frequency interval allows a Taylor expansion of the cosine in \eqref{eq:DFourier}. Its constant term gives the density at zero. Its quadratic term will decrease the density provided that its coefficient has the right sign after averaging over \(B\). For \(B\in\Hcal\), write
\begin{equation}\label{eq:AH}
\begin{gathered}
 D_B(0)=2\int_{\R^m}\frac{p_B^{\rm pol}(0,x)}{W_0^*(x)}\,dx,\qquad
 \mathcal K_B=\int_{-K}^K u^2\operatorname{Re}\E[e^{iuL}W_0^*(Z)\mid B]\,du,\\
 \bar D_n=\E[\one{B\in\Hcal} D_B(0)] .
\end{gathered}
\end{equation}
In expectations restricted to a block event, define the integrand to be zero off that event. Here \(D_B(0)\) uses the smoothing and the inverse weight of this section. The coefficient \(C_B\) in \eqref{eq:flatcoef} used a different weight. The density of \(L\) at zero is \(O(1)\) uniformly on \(\Hcal\). Indeed Fourier inversion gives an integrand bounded by \(|\varphi(uv)|e^{-n^{-180}u^2/2}\). Integrate Corollary~\ref{co:slab} up to \(\sqrt n\), use \(e^{-c_1 n}\) between \(\sqrt n\) and \(n^{91}\), and integrate the Gaussian tail thereafter. Thus \(0\le D_B(0)\le C\sqrt n\). Also
\begin{equation}\label{eq:Wmoms}
 \E[W_0^*(Z)\mid B]\le C\sqrt{\log n}(1+a_1^{-1}),\qquad
 \E[W_0^*(Z)^3\mid B]\le C\log^{3/2}n(1+a_1^{-3}).
\end{equation}
For example, apply Minkowski in \(L^{3/2}\) to \(b+\sum Z_j^2\) and use \(\|Z_j\|_{L^3}\le C(a_j^{-1}+r)\), \eqref{eq:tails}, and \(mr^2=o(1)\).

\begin{lemma}\label{le:low}
For \(B\in\Hcal_t,\ t_*\le t\le n^{-\theta},\ 0\le s\le2t\),
\begin{align}
 D_B(s)&\le D_B(0)-\frac{s^2}{2\pi}\mathcal K_B
       +Cs^4K^5\log^{3/2}n(1+a_1^{-3})+3n^{-100},\label{eq:Dupper}\\
 D_B(s)&\ge (1-s^2/a_1^2)D_B(0)
       -Cs^2K^3\sqrt{\log n}(1+a_1^{-1})-2n^{-100}.\label{eq:Dlower}
\end{align}
\end{lemma}
\begin{proof}
Suppress the argument \(Z\) in \(W_s^*(Z)\) and \(W_0^*(Z)\). Then \(W_0^*\le W_s^*\le(1+s^2/a_1^2)W_0^*\le1.01W_0^*\). By Fourier inversion as in Lemma~\ref{le:crossdensity},
\[
 \frac{1}{\pi}\int_\R\operatorname{Re}\E[e^{iuL}/W_s^*\mid B]\,du
=2\int_{\mathbb R^m}\frac{p_B^{\rm pol}(0,x)}{W_s^*(x)}\,dx.
\]
By positivity this integral lies between \((1-s^2/a_1^2)D_B(0)\) and \(D_B(0)\). Lemma~\ref{le:high} permits truncation of this integral and \eqref{eq:DFourier} at \(K\). For \(|u|\le K\), \(|\cos(usW_s^*)-1|\le u^2s^2(W_s^*)^2/2\), proving \eqref{eq:Dlower} by \eqref{eq:Wmoms}. For the upper bound, expand \(\cos(usW_s^*)=1-u^2s^2(W_s^*)^2/2+\mathcal E_{\cos}\), \(|\mathcal E_{\cos}|\le u^4s^4(W_s^*)^4/24\). The integrated remainder is at most \(Cs^4K^5\E[(W_s^*)^3\mid B]\). Replacing \(W_s^*\) by \(W_0^*\) in the quadratic term costs at most \(Cs^4 K^3a_1^{-2}\E[W_0^*\mid B]\). Both errors are absorbed in \eqref{eq:Dupper} by \eqref{eq:Wmoms}.
\end{proof}
\begin{lemma}[Curvature]\label{le:curv}
On \(B\in\Hcal\),
\begin{enumerate}[label=\textup{(\alph*)}]
\item \(|\mathcal K_B|\le C K^3\sqrt{\log n}(1+a_1^{-1})\).
\item If \(a_1\ge n^{-\theta/2}\), then \(\mathcal K_B\ge c_{\rm curv}/\log n\), where \(c_{\rm curv}=1/(2C_*)\), for sufficiently large \(n\).
\end{enumerate}
\end{lemma}
\begin{proof}
Part (a) follows by taking absolute values in \eqref{eq:AH} and applying \eqref{eq:Wmoms}. For (b) we need the sign of the integral. The corresponding Gaussian calculation separates into two positive factors because the normal and tangential projections are independent. We show that replacing the sign row by a Gaussian row changes this integral by \(o(1/\log n)\). Fix a block with \(a_1\ge n^{-\theta/2}\).

\step{Reduce to the first \(k\) tangential coordinates}
On \([-K,K]\) removing the Gaussian factor \(e^{-n^{-180}u^2/2}\) changes \(\mathcal K_B\) by at most \(Cn^{-180}K^5\E W_0^*(Z)\le n^{-170}\).
Define
\(\rho(\xi')=(n^{-2\kappa}+b+|\xi'|^2)^{1/2}\), \(\xi'\in\R^k\), and let \(V=\rho(z_1,\ldots,z_k)\). The triangle inequality and \eqref{eq:tails} yield
\[
 \E|W_0^*(Z)-V|\le n^{-\kappa}
       +\Big(\sum_{k<j\le m}a_j^{-2}\Big)^{1/2}+C r\sqrt m
       \le C n^{-\kappa/2}.
\]
The resulting error in \(\mathcal K_B\) is at most \(C K^3 n^{-\kappa/2}\).

\step{Replace one sign entry at a time} For fixed \(|u|\le K\), set \(\Phi(\xi_0,\xi')=e^{i\xi_0}\rho(\xi')\) and \(h_\iota=(uv_\iota,(u_{j,\iota}/a_j)_{j\le k})\). Then \(e^{iu\ell}V=\Phi(\sum_\iota X_\iota h_\iota)\). The first three derivatives of \(\rho\) have operator norms at most \(1,n^\kappa,6n^{2\kappa}\), respectively. Hence
\[
 |D^3\Phi(\xi)[h,h,h]|\le(\rho(\xi')+12n^{2\kappa})|h|^3,\qquad
 \rho(\xi'+h')\le\rho(\xi')+|h'|.
\]
Replace \(X_\iota\) successively by independent standard Gaussians \(g_\iota\). In a single replacement let \(Y_{-\iota}\) be the hybrid sum excluding index \(\iota\), so that \(Y_{-\iota}\) is independent of both scalars being compared. Taylor's formula, cancellation of the first two moments and finiteness of the third and fourth absolute moments give a total error of at most
\[
 C\sum_\iota |h_\iota|^3
       \big(\E\rho(Y_{-\iota}')+n^{2\kappa}+|h_\iota|\big).
\]
Put \(\mathcal E_h=\sum_\iota|h_\iota|^2=u^2+\sum_{j\le k}a_j^{-2}\le C n^\theta\log n\). One has \(\E\rho(Y_{-\iota}')\le C(n^\theta\log n)^{1/2}\) by the second moments. By (G5), each coordinate of \(v\) and of each \(u_j\) for \(j\le k\) is at most \(n^{-1/4}\) in absolute value, so \(\max_\iota |h_\iota|^2\le n^{-1/2}\mathcal E_h\). Thus the total error is \(O(n^{-1/4}(n^\theta\log n)^2)\le n^{-1/5}\), uniformly over \(|u|\le K\).

\step{Evaluate the Gaussian integral} For Gaussian input, the projections onto \(v,u_1,\ldots,u_k\) are independent. Consequently
\[
 \E\Phi\Big(\sum_\iota g_\iota h_\iota\Big)=e^{-u^2/2}\E\rho(\gamma_1/a_1,\ldots,\gamma_k/a_k)
\]
for iid standard Gaussians \(\gamma_j\). The square root in \(\rho\) is at least the absolute first coordinate, so the second factor is at least \(\sqrt{2/\pi}/a_1\ge\sqrt{2/\pi}/(C_*\log n)\) by (G3). For large \(n\), \(\int_{-K}^K u^2e^{-u^2/2}du\ge\sqrt{2\pi}/2\), so
\[
 \mathcal K_B\ge \frac1{C_*\log n}
    -CK^3(n^{-\kappa/2}+n^{-1/5})-n^{-170}
 \ge \frac1{2C_*\log n}
\]
since \(3\tau<\kappa/2\).
\end{proof}

\subsection{Proof of Proposition~\ref{pr:poly}}
We first average the upper density bound over the retained blocks. Its negative quadratic term will pay for the discarded blocks when we integrate over the threshold. A separate lower bound at \(t_*\) then determines how large the averaged density at zero can be.

\step{Average the curvature} Write \(T'=\max(T,n^{-4})\le n^{-\theta/2}\) for large \(n\). In \eqref{eq:Dupper}, positivity gives \(\E[\one{B\in\Hcal_t}D_B(0)]\le \bar D_n\). The good blocks with \(a_1\ge n^{-\theta/2}\) all belong to \(\Hcal_t\). Their total probability is at least \(1-n^{-100}-C n^{-\theta}\) by Proposition~\ref{pr:good} and Lemma~\ref{le:RV}(a). For \(x_0=n^{-\theta/2}\), Lemma~\ref{le:RV}(a) gives
\[
 \E[a_1^{-1}\one{n^{-4}\le a_1<x_0}]
 \le x_0^{-1}\PP(a_1\le x_0)
       +\int_{n^{-4}}^{x_0}\PP(a_1\le x)x^{-2}\,dx
 \le Cx_0 .
\]
Lemma~\ref{le:curv} and \(3\tau<\theta/2\) now give
\(\E[\one{B\in\Hcal_t}\mathcal K_B]\ge c_{\rm curv}/(2\log n)\).
\step{Average the Taylor remainder} For the quartic remainder in \eqref{eq:Dupper}, integrate Lemma~\ref{le:RV}(a) against \(3x^{-4}\). This gives
\[
 \E[\one{B\in\Hcal_t}(1+a_1^{-3})]
 \le 1+\int_{T'}^\infty 3x^{-4}\PP(a_1\le x)\,dx
 \le C/T'\le C/(t\log^3 n),
\]
using \(T'\ge n^{-4}\) to absorb the exponential error. For \(s\le2t\) this bounds the averaged remainder by \(Cs^4K^5/(t\log^{3/2}n)\), which is at most \(Cs^2 n^{5\tau-\theta}\). With \(c_{\rm loss}=c_{\rm curv}/(8\pi)\), this proves
\begin{equation}\label{eq:neg}
 \E[\one{B\in\Hcal_t}D_B(s)]\le \bar D_n-c_{\rm loss}s^2/\log n+3n^{-100},
 \qquad 0\le s\le2t .
\end{equation}
\step{Integrate the density and restore the discarded blocks} Combine Lemmas~\ref{le:exchange}, \ref{le:smallblock}, \ref{le:compare}, and \ref{le:crossdensity}. Integrating \eqref{eq:neg} up to \((1+\epsilon_n)t\) gives
\[
 F_n(t)\le (1+\epsilon_n)^2 \bar D_n t
      -c_{\rm loss}t^3/(3\log n)+t^{10/3}n^{3\kappa}+2tn^{-7}.
\]
Indeed \(n^{-90}\) and the other additive errors are absorbed since \(t\ge n^{-20}\). The ratio of \(t^{10/3}n^{3\kappa}\) to \(t^3/\log n\) is at most \(n^{-\theta/3+3\kappa}\log n=o(1)\). The two terms together are nonpositive for sufficiently large \(n\). Since \((1+\epsilon_n)^2\le1+3\epsilon_n\), this proves (a).

\step{The anchor \(t_*\)} For \(B\in\Hcal_{t_*}=\Hcal\) and \(s\le2t_*\), \(a_1\ge n^{-4}\) and \(D_B(0)\le C\sqrt n\). Substitution in \eqref{eq:Dlower} gives \(D_B(s)\ge D_B(0)-n^{-30}\). Integrate to \((1-\epsilon_n)t_*\), use the lower bound of Lemma~\ref{le:compare} on \(\Hcal\), and discard the nonnegative contribution of the other blocks. Then
\[
 F_n(t_*)\ge(1-\epsilon_n)^2t_*\bar D_n-t_*n^{-30}-n^{-90}
       \ge(1-3\epsilon_n)t_* \bar D_n-t_*n^{-6}.
\]
This proves (b).
\section{Proof of the theorems}\label{sec:finish}
\subsection{Choice of the polynomial threshold}
Tao and Vu prove a uniform universality estimate \cite[Theorem~1.3, eq.~(8) in the arXiv version]{taodist}. Applied to real Rademacher entries, it gives
\begin{equation}\label{eq:TV}
 |F_n(t)-F_\infty(t)|\le C_{\rm TV}n^{-c_{\rm TV}},\qquad
 F_\infty(t)=1-e^{-t-t^2/2}.
\end{equation}
Their variable is \(n\smin(A)^2\). Integrating their limiting density up to \(t^2\) and substituting \(x=y^2\) gives \(\int_0^t(1+y)e^{-y-y^2/2}\,dy=F_\infty(t)\). The error is uniform in the threshold, so this estimate also holds when \(t\) depends on \(n\). The result for positive thresholds extends to zero by taking a decreasing limit. Choose \(\theta=\min(c_{\rm TV}/4,1/100)\).

\subsection{Proof of Theorem~\ref{th:sharper}}
Set \(c=\cflat\). On \([e^{-\cflat n},n^{-20}]\), Proposition~\ref{pr:flat} and Corollary~\ref{co:slope} give \(F_n(t)\le t(1-1/(2n)+O(n^{-1.05}))\). On \([n^{-20},n^{-\theta}]\) use Corollary~\ref{co:poly}.

\subsection{Proof of Theorem~\ref{th:main}}
\step{Large dimensions} For \(n\) large enough, Theorem~\ref{th:sharper} yields \(F_n(t)\le t\) on \([e^{-\cflat n},n^{-\theta}]\), because \(Cn^{-1.02}<1/(2n)\). For smaller \(t\ge0\), monotonicity gives \(F_n(t)\le F_n(e^{-\cflat n})\le e^{-\cflat n}\). On \([0,1]\) the function \(g(t)=t-F_\infty(t)\) has \(g(0)=g'(0)=0\) and \(g''(t)=(2t+t^2)e^{-t-t^2/2}\ge2e^{-3/2}t\). Hence, for \(n^{-\theta}\le t\le1\),
\[
 F_n(t)\le t-\tfrac13e^{-3/2}n^{-3\theta}+C_{\rm TV}n^{-4\theta}\le t
\]
for large \(n\). For \(t\ge1\) use \(F_n(t)\le1\).

\step{The remaining dimensions} Fix an integer \(N\) beyond which the preceding argument holds. To obtain one exponential rate for all \(n\le N\), set \(\eta_n=2^{-n^2}\), the probability of a prescribed sign matrix. An invertible sign matrix exists in each dimension. Indeed the all-ones matrix minus \(2I\) is invertible for \(n\ne2\), and for \(n=2\) take rows \((1,1)\) and \((1,-1)\). Hence \(\PP(X_n>0)\ge\eta_n\).

In an invertible realization the determinant is a nonzero integer and every singular value is at most \(\|A\|_F=n\). From their product we obtain \(X_n\ge\sqrt n\,n^{-(n-1)}\ge\eta_n\). For \(t<\eta_n\) it follows that \(F_n(t)=F_n(0)\le1-\eta_n\). For \(t\ge\eta_n\) use \(F_n(t)-t\le1-t\le1-\eta_n\). Taking
\[
 c=\min\{\cflat,\, 2^{-N^2}/N\}
\]
now works in every dimension. If \(n\le N\), then \(cn\le\eta_n\) and \(1-\eta_n\le e^{-cn}\). For \(n>N\), decreasing the rate only enlarges the permitted exponential error.

\begin{remark}
This formula makes the finite-dimensional step explicit once \(N\) and \(\cflat\) are given. The argument does not calculate either value. For example, Lemma~\ref{le:high} uses \(\log n/n^{\tau/2}\to0\) with \(\tau=\theta/1000\). Lemma~\ref{le:many} uses convergence of an expected eigenvalue count without a numerical rate.
\end{remark}
\begin{remark}
For independent, centered, variance-one entries, the almost sure row-length bound \(|X|\le\sqrt n\) forces each coordinate to be a uniform sign. Indeed \(\E|X|^2=n\) forces the squared length to equal \(n\) almost surely, and independence then forces each \(X_i^2\) to be constant. Lemma~\ref{le:angle} also uses the sign characteristic function. The coefficient \(1/2\) in Theorem~\ref{th:sharper} comes from separate upper bounds on two factors. Corollary~\ref{co:slope} gives \(\liminf_{n\to\infty} n(1-a_n)\ge1/2\). The existence and value of a limit remain undetermined by this argument.
\end{remark}
\appendix
\section{The spherical ensemble}\label{app:sph}
Assume \(n\ge8\). Each of the \(n\) rows of \(A_{\sph}\) is sampled independently and uniformly from \(\sqrt n\,\sphere\). Let \(F_{\sph}\) denote the distribution function of \(\sqrt n\,\smin(A_{\sph})\).
\begin{proposition}\label{pr:sph}
The right derivative at zero exists and satisfies
\[
 F_{\sph}'(0+)=\frac{|\sphere|}{\sqrt n}\,f_n^{\,n}\E J_*
     =1-\frac{5}{4n}+O(n^{-2}).
\]
\end{proposition}

The two equalities in Proposition~\ref{pr:sph} serve different purposes. To find the numerical coefficient, we first recover the spherical distribution from a Gaussian matrix by conditioning on its row lengths. We then differentiate its survival probability. To identify the coefficient with the Jacobian in Section~\ref{sec:area}, we compute the same derivative by exposing the last row.

\subsection{A Gaussian identity}
\begin{lemma}\label{le:wish}
Let \(\Gamma\) have iid standard Gaussian entries, with \(n\) rows and \(n\) columns, let \(z\sim N(0,I_n)\), and let \(\Xi\) be deterministic and positive semidefinite. Then
\[
 \PP(\Gamma\Gamma^\top\succ \Xi)=e^{-\tr \Xi/2}\E e^{-|\Xi^{1/2}z|}.
\]
Here \(\succ\) means that the difference is positive definite.
\end{lemma}
\begin{proof}
The density of \(\mathsf W=\Gamma\Gamma^\top\) on the positive definite cone is proportional to \(\det(\mathsf W)^{-1/2}e^{-\tr \mathsf W/2}\) \cite[Theorem~3.2.1, p.~85]{muirhead}. Shifting by \(\Xi\) in that integral gives
\[
 \PP(\mathsf W\succ \Xi)=e^{-\tr \Xi/2}\E\det(I+\Xi^{1/2}\mathsf W^{-1}\Xi^{1/2})^{-1/2}.
\]
Use \(\det(I+M)^{-1/2}=\E_z e^{-z^\top Mz/2}\) for \(M\succeq0\), taking \(z\) independent of \(\mathsf W\). For fixed \(y\ne0\), \(y^\top \mathsf W^{-1}y\) has law \(|y|^2/g^2\) for standard normal \(g\). This follows by rotating \(y\) to the last coordinate and using the Schur complement of the first \(n-1\) rows. The reciprocal of the last diagonal entry of \(\mathsf W^{-1}\) is the squared distance of the last Gaussian row to their span, which has law \(\chi^2_1\). Finally \(\E e^{-\alpha^2/(2g^2)}=e^{-\alpha}\) for \(\alpha\ge0\), by the Gaussian integral \(\int_0^\infty e^{-p x^2-q/x^2}\,dx=\frac12\sqrt{\pi/p}\,e^{-2\sqrt{pq}}\) for \(p>0,\ q\ge0\). Apply this with \(y=\Xi^{1/2}z\), including \(y=0\) by continuity.
\end{proof}
\begin{remark}\label{rm:Gaussian}
Taking \(\Xi=(t^2/n)I\) yields \(\PP(\sqrt n\,\smin(\Gamma)\le t)=1-e^{-t^2/2}\E e^{-t|z|/\sqrt n}\le1-e^{-t-t^2/2}\le t\). The first inequality follows from Jensen's inequality and \(\E|z|\le\sqrt n\). The second follows from Section~\ref{sec:finish}.
\end{remark}

\subsection{The exact spherical law}
\begin{lemma}\label{le:exactsph}
Let \(z_1,\ldots,z_n,g\) be independent standard Gaussians. For \(0<t<n\), put \(\eta_t=t^2/n^2\). Then
\[
 1-F_{\sph}(t)=\E\prod_{i=1}^n
      \left(1-\eta_t\left(1+\frac{z_i^2}{g^2}\right)\right)_+^{n/2-1}.
\]
The same \(g\) appears in every factor.
\end{lemma}
\begin{proof}
A Gaussian row has a uniform direction independent of its length. We first allow each squared row length to vary, so the Gaussian identity tests the unknown spherical law against a family of gamma densities.

Let \(\Omega\) have independent rows uniform on \(\sphere\). Fix \(s>0\). For \(h_i>0\) define
\[
 \mathcal S_s(h)=\PP(\smin(\diag(\sqrt{h_i})\Omega)>s),\qquad
 \widetilde{\mathcal S}_s(h)=\E\prod_i\left(1-\frac{s^2}{h_i}\left(1+\frac{z_i^2}{g^2}\right)\right)_+^{n/2-1}.
\]
For any \(\lambda_i>0\), the matrix with rows \(\Gamma_i/\sqrt{2\lambda_i}\) has the law of \(\diag(\sqrt{h_i})\Omega\) with independent \(h_i\sim\operatorname{Gamma}(n/2,\text{rate }\lambda_i)\), independent of \(\Omega\). Its least singular value exceeds \(s\) exactly when \(\Gamma\Gamma^\top\succ 2s^2\diag(\lambda_i)\). Lemma~\ref{le:wish} and the scalar Gaussian integral used there give
\[
 \E_h \mathcal S_s(h)
 =e^{-s^2\sum_i\lambda_i}\E_z\exp\!\left(-s\sqrt{2\sum_i\lambda_i z_i^2}\right)
 =\E_{z,g}\prod_i e^{-\lambda_i s^2(1+z_i^2/g^2)} .
\]
For \(h\) gamma of shape \(b=n/2\) and rate \(\lambda\), multiplication by \((1-\beta/h)_+^{b-1}\) simply shifts the power in its density. For \(\beta\ge0\),
\[
 \E(1-\beta/h)_+^{b-1}
 =\frac{\lambda^b}{\Gamma(b)}\int_\beta^\infty
       (h-\beta)^{b-1}e^{-\lambda h}\,dh=e^{-\lambda\beta}.
\] Therefore \(\E_h\widetilde{\mathcal S}_s(h)\) equals the same expression for every positive \(\lambda\). The functions \(\mathcal S_s(h)\prod_i h_i^{n/2-1}\) and \(\widetilde{\mathcal S}_s(h)\prod_i h_i^{n/2-1}\) are locally integrable and integrable after every strictly positive exponential tilt. Uniqueness of the multivariate Laplace transform (or successive one-dimensional uniqueness) yields \(\mathcal S_s=\widetilde{\mathcal S}_s\) almost everywhere.

For each realization, \(\Omega^\top\diag(h_i)\Omega\) increases in the positive semidefinite order when any \(h_i\) increases. Hence \(\mathcal S_s\) is coordinatewise nondecreasing. The function \(\widetilde{\mathcal S}_s\) is continuous by dominated convergence. For any \(h_0\in(0,\infty)^n\), choose sequences in the full-measure equality set approaching \(h_0\) from strictly below and strictly above in every coordinate. Monotonicity squeezes \(\mathcal S_s(h_0)\) between values converging to \(\widetilde{\mathcal S}_s(h_0)\). Set \(h_0=(n,\ldots,n)\) and \(s=t/\sqrt n\) to obtain the formula.
\end{proof}

\begin{lemma}\label{le:sphder}
\(F_{\sph}'(0+)=1-5/(4n)+O(n^{-2})\).
\end{lemma}
\begin{proof}
\step{Differentiate at fixed dimension} Write the product in Lemma~\ref{le:exactsph} as a probability by introducing one auxiliary variable per row. Set \(m_0=n/2-1\) and take iid \(V_i\) independent of \((z,g)\) with \(\PP(V_i>x)=(1-x)_+^{m_0}\) for \(x\ge0\). Lemma~\ref{le:exactsph} expresses the survival probability as \(\PP(V_i>\eta_t+\eta_t z_i^2/g^2\ \forall i)\). Given \(V_i>\eta_t\) for all \(i\), the rescaled variables \(V_i'=(V_i-\eta_t)/(1-\eta_t)\) are again iid with the original law. Hence, using fresh copies in the following expression,
\[
 1-F_{\sph}(t)=(1-\eta_t)^{nm_0}\PP\!\left(|g|>\sqrt{\frac{\eta_t}{1-\eta_t}}M_n\right),
 \qquad M_n=\max_{i\le n}|z_i|/\sqrt{V_i'} .
\]
Since \(\PP(V_i'\le x)\le m_0x\) for \(x\ge0\), conditioning on each \(z_i\) and taking a union bound gives \(\PP(M_n>ny)\le nm_0/(n^2y^2)\le1/(2y^2)\) for \(y>0\). In particular \(\E M_n<\infty\). At fixed \(n\), the prefactor is \(1-O(t^2)\) and \(\PP(|g|\le h M_n)/h\to \sqrt{2/\pi}\,\E M_n\) as \(h\downarrow0\) by dominated convergence. The spherical matrix has full rank almost surely, hence \(F_{\sph}(0)=0\). It follows that
\[
 F_{\sph}'(0+)=\sqrt{\frac2\pi}\frac{\E M_n}{n}
 =\sqrt{\frac2\pi}\int_0^\infty\frac{1-\psi_n(x)^n}{x^2}\,dx,\qquad
 \psi_n(x)=\E(1-x^2(g')^2/n^2)_+^{m_0}.
\]
Here \(g'\) denotes a standard Gaussian. The integral follows from \(\PP(M_n\le ny)=\psi_n(1/y)^n\).

\step{Extract the first finite-dimensional correction} We now let \(n\) tend to infinity in the integral for the derivative. On \(y\ge0\), the second derivative of \((1-y)_+^{m_0}\) is Lipschitz with constant at most \(m_0^3\). Taylor expansion and the Gaussian moments \(1,3,15\) give
\[
 \psi_n(x)=1-\frac{m_0x^2}{n^2}+\frac{3m_0(m_0-1)x^4}{2n^4}
       +O(x^6/n^3)
 =1-\frac{x^2}{2n}+\frac{x^2+3x^4/8}{n^2}
       +O((x^4+x^6)/n^3).
\]
For \(0\le x\le x_1=\sqrt{20\log n}\) the logarithm can be expanded uniformly. The quadratic term of \(\log(1+u)\) subtracts \(x^4/(8n)\) from \(n(\psi_n(x)-1)\), changing the coefficient of \(x^4/n\) from \(3/8\) to \(1/4\). Exponentiation gives
\[
 \psi_n(x)^n=e^{-x^2/2}\left[1+\frac{x^2+x^4/4}{n}
                   +O\!\left(\frac{x^4+x^8}{n^2}\right)\right].
\]
The error after division by \(x^2\) is integrable uniformly, including at zero. At the right endpoint \(\psi_n(x_1)^n\le n^{-9}\) for large \(n\). To control the rest of the integral, subtract the Gaussian baseline first. Monotonicity of \(\psi_n\) gives
\[
 \int_{x_1}^\infty \frac{|\psi_n(x)^n-e^{-x^2/2}|}{x^2}\,dx
 \le \frac{n^{-9}+e^{-x_1^2/2}}{x_1}.
\]
Indeed, the two integrands \((1-\psi_n(x)^n)/x^2\) and \((1-e^{-x^2/2})/x^2\) have the same \(x^{-2}\) tail. Integration by parts shows \(\int_0^\infty(1-e^{-x^2/2})x^{-2}\,dx=\sqrt{\pi/2}\). Thus
\[
 \int_0^\infty\frac{1-\psi_n(x)^n}{x^2}\,dx
 =\sqrt{\pi/2}-\frac1n\int_0^\infty e^{-x^2/2}(1+x^2/4)\,dx+O(n^{-2})
 =\sqrt{\pi/2}\left(1-\frac5{4n}+O(n^{-2})\right).
\]
This proves the expansion.
\end{proof}

\subsection{Identification of the geometric constant}
\begin{lemma}\label{le:sphgeom}
\(F_{\sph}'(0+)=|\sphere|\,n^{-1/2} f_n^{\,n}\E J_*\).
\end{lemma}
\begin{proof}
\step{Conditional derivative} Split \(A_{\sph}\) into \(B_s\) and its last row \(X_s\) as in Lemma~\ref{le:schur}. The block has full row rank almost surely. Conditionally on \(B_s\), the normal coordinate \(\ell\) of \(X_s\) has density \(f(l)=f_n(1-l^2/n)_+^{(n-3)/2}\). Given \(\ell\), the vector \((y_j)\) has law \(\sqrt{n-\ell^2}\,U'\) for \(U'\) uniform on \(S^{n-2}\) and independent of \(\ell\). For \(0<h<d_1\) put \(\mathcal R_h=\sum_j (U'_j)^2/(d_j^2-h^2)\). By Lemma~\ref{le:schur}(ii) the event \(\smin(A_{\sph})\le h\) is equivalent to \(\ell^2\le h^2[1+(n-\ell^2)\mathcal R_h]\). Gathering the two terms involving \(\ell^2\), we obtain \(|\ell|\le h\omega_h\), where \(\omega_h=((1+n\mathcal R_h)/(1+h^2 \mathcal R_h))^{1/2}\). For \(h\le d_1/2\), \(\omega_h\le1+2\sqrt n/d_1\). Hence
\begin{equation}\label{eq:condder}
 \lim_{h\downarrow0}\frac{\PP(\smin(A_{\sph})\le h\mid B_s)}{h}
   =2f_n\Phi(B_s),\qquad
 \Phi(B_s)=\E_{U'}\Big(1+n\sum_j(U'_j)^2/d_j^2\Big)^{1/2}.
\end{equation}
For \(h\le d_1/2\), the conditional quotient is bounded by \(2f_n(1+2\sqrt n/d_1)\). For larger \(h\) it is at most \(2/d_1\). To justify averaging at each fixed \(n\), let \(\Gamma\) be \((n-1)\times n\) standard Gaussian. Under row normalization \(d_1\ge\sqrt n\,s_1(\Gamma)/\max_i|\Gamma_i|\). For a full-rank rectangular matrix the smallest singular value is at least the least row-to-other-rows distance divided by \(\sqrt{n-1}\). To see this, take a unit combination of rows and project off all but a row with coefficient of largest magnitude. Each squared distance for \(\Gamma\) has law \(\chi_2^2\). Thus \(\PP(s_1(\Gamma)\le x)\le C n^2x^2\) and \(\E s_1(\Gamma)^{-4/3}<\infty\). Hölder's inequality gives \(\E d_1^{-1}<\infty\). Dominated convergence now yields
\begin{equation}\label{eq:unconder}
 \lim_{h\downarrow0}\frac{\PP(\smin(A_{\sph})\le h)}h=2 f_n\E\Phi(B_s).
\end{equation}

\step{Compare the two row dimensions} The rows of \(B_s\) live in \(\R^n\). The rows of the Jacobian matrix \(Y_*\) live in \(\R^{n-1}\). The following identity compares their Gram laws after scaling to unit row lengths. Let \(\Omega'=B_s/\sqrt n\). Let \(\Theta\) be an \((n-1)\times(n-1)\) matrix with iid rows uniform on \(S^{n-2}\). For every nonnegative measurable function \(\chi\) of a Gram matrix,
\begin{equation}\label{eq:bias}
 \E\chi(\Omega'\Omega'^\top)
   =\frac{\E[|\det\Theta|\,\chi(\Theta\Theta^\top)]}{\E|\det\Theta|}.
\end{equation}
Indeed the Gram densities for Gaussian matrices of sizes \((n-1)\times n\) and \((n-1)\times(n-1)\) are proportional to \(e^{-\tr \mathsf W/2}\) and \(\det(\mathsf W)^{-1/2}e^{-\tr \mathsf W/2}\), respectively. Thus the first is the second tilted by \(\sqrt{\det \mathsf W}\). Apply this to the correlation matrix \(\diag(\mathsf W)^{-1/2}\mathsf W\diag(\mathsf W)^{-1/2}\), where \(\diag(\mathsf W)\) denotes the diagonal matrix with the diagonal of \(\mathsf W\). For the square Gaussian matrix the tilting factor is \(|\det\Theta|\) times the product of the independent row lengths. Those lengths are independent of the directions and cancel upon normalization.

\step{Cancel the determinant weight} Since \(\Phi\) depends only on the singular values, \eqref{eq:bias} applies. For an invertible square \(M\) of size \(n-1\), the determinant lemma and rotational invariance of \(U'\) give
\[
 |\det M|\,\E_{U'}\Big(1+n\sum_j (U_j')^2/s_j(M)^2\Big)^{1/2}
 =\E_{U'}\sqrt{\det(M^\top M+nU'(U')^\top)} .
\]
With \(M=\sqrt n\,\Theta\), appending \(\sqrt n\,(U')^\top\) gives the law of \(Y_*\). Hence
\(\E\Phi(B_s)=\E J_*/(n^{(n-1)/2}\E|\det\Theta|)\).
For an \((n-1)\times(n-1)\) Gaussian matrix \(\Gamma'\), Gram--Schmidt gives \(\E|\det\Gamma'|=\prod_{i=1}^{n-1}\E\chi_i\), while independence of the row radii gives \(\E|\det\Gamma'|=(\E\chi_{n-1})^{n-1}\E|\det\Theta|\), where \(\E\chi_i=\sqrt2\,\Gamma((i+1)/2)/\Gamma(i/2)\). Simplifying,
\[
 \E|\det\Theta|=\frac{\Gamma(n/2)}{\sqrt\pi}
       \left(\frac{\Gamma((n-1)/2)}{\Gamma(n/2)}\right)^{n-1},
 \qquad n^{(n-1)/2}\E|\det\Theta|=\frac{2}{|\sphere|f_n^{\,n-1}} .
\]
Substitute in \eqref{eq:unconder} and use \(F_{\sph}(t)=\PP(\smin(A_{\sph})\le t/\sqrt n)\). This proves the identity and, with Lemma~\ref{le:sphder}, Proposition~\ref{pr:sph}.
\end{proof}

\end{document}